\documentclass[a4paper,reqno,11pt]{amsart}
\usepackage[english]{babel}
\usepackage{amsmath}
\usepackage{amssymb}
 \usepackage{amsthm}
\usepackage{enumerate}
\usepackage{ifthen}
\usepackage{bbm}
\usepackage{color}
\usepackage{cancel}
\usepackage{hyperref}
\usepackage{geometry}
\newtheorem{theorem}{Theorem}[section]
\newtheorem{proposition}[theorem]{Proposition}
\newtheorem{lemma}[theorem]{Lemma}

\newtheorem{definition}[theorem]{Definition}
\theoremstyle{remark}
\newtheorem{remark}[theorem]{Remark}

\newcommand{\R}{\mathbb{R}}
\newcommand{\tore}{\mathbb{T}^3}
\newcommand{\N}{\mathbb{N}}
\newcommand{\Z}{\mathbb{Z}}
\newcommand{\dive}{\mathop{\mathrm {div}}}

\newcommand{\weakto}{\rightharpoonup}
\newcommand{\pth}{p_h^{m}}
\newcommand{\uth}{u_h^{m,1/2}}

\numberwithin{equation}{section}
\subjclass[2010]{Primary: 35Q30, Secondary: 
65M12, 
    76M20.}
\keywords{Navier-Stokes Equations, Local energy inequality, Numerical
  schemes, second-order methods, Finite Element and Finite Difference
  Methods}
\begin{document}
\title[Convergence of second order schemes for the 
NSE]{On the convergence of second-order in time numerical
  discretizations for the evolution Navier-Stokes equations}
\author[Luigi C. Berselli]{Luigi C. Berselli} 
\address[L. C. Berselli] {Dipartimento di Matematica
  \\
  Universit\`a degli Studi di Pisa, Via F. Buonarroti 1/c, I-56127 Pisa, Italy}
\email[]{\href{luigi.carlo.berselli@unipi.it}{luigi.carlo.berselli@unipi.it}}
\author[Stefano. Spirito]{Stefano Spirito} 
\address[S. Spirito]{DISIM - Dipartimento di Ingegneria e Scienze
  dell'Informazione e Matematica
  \\
  Universit\`a degli Studi dell'Aquila, Via Vetoio, I-67100 L'Aquila,
  Italy.}
\email[]{\href{stefano.spirito@univaq.it}{stefano.spirito@univaq.it}}

\begin{abstract}
  We prove the convergence of certain second-order numerical methods
  to weak solutions of the Navier-Stokes equations satisfying in
  addition the local energy inequality, and therefore suitable in the
  sense of Scheffer and Caffarelli-Kohn-Nirenberg. More precisely, we
  treat the space-periodic case in three space-dimensions and we
  consider a full discretization in which the the classical
  Crank-Nicolson method ($\theta$-method with $\theta=1/2$) is used to
  discretize the time variable, while in the space variables we
  consider finite elements. The convective term is discretized in
  several implicit, semi-implicit, and explicit ways.  In particular,
  we focus on proving (possibly conditional) convergence of the
  discrete solutions towards weak solutions (satisfying a precise
  local energy balance), without extra regularity assumptions on the
  limit problem. We do not prove orders of convergence, but our
  analysis identifies some numerical schemes providing also alternate
  proofs of existence of ``physically relevant'' solutions in three
  space dimensions.
\end{abstract}
\maketitle
\section{Introduction}
We consider the homogeneous incompressible 3D~Navier-Stokes equations
(NSE)
\begin{equation}
  \label{eq:nse}
  \begin{aligned}
    \partial_t u-\nu\Delta u+(u\cdot\nabla)\,u+\nabla
    p&=0\qquad\textrm{ in }     (0,T)\times\tore,
    \\
    \dive u&=0\qquad\textrm{ in }
    (0,T)\times\tore,
     \end{aligned}
\end{equation}
in the space periodic setting, with divergence-free initial datum
 \begin{equation}
   \label{eq:nsid}
 u|_{t=0}=u_0\qquad\textrm{ in }\tore,
 \end{equation}
 where $T>0$ is arbitrary, $\nu>0$ is given, and
 $\tore :=(\R/ 2\pi\Z)^3$ is the three dimensional flat torus. Here,
 the unknowns are the velocity vector field $u$ and the scalar
 pressure $p$, which are both with zero mean value.  The aim of this
 paper is to consider families of space-time discretization of the
 initial value problem~\eqref{eq:nse}-\eqref{eq:nsid} which are of
 second order in time and (as the parameters of the discretization
 vanish) to prove the convergence towards Leray-Hopf weak solutions,
 satisfying in addition certain estimates on the pressure and the
 local energy inequality
 \begin{equation*}
   \partial_t\left(\frac{|u|^2}{2}\right)+\dive\left(\left(\frac{|u|^2}{2}
       +p\right)u\right)-   \nu\Delta
   \left(\frac{|u|^2}{2}\right)+\nu|\nabla u|^2\leq 0\qquad 
   \text{in }\mathcal{D}'(]0,T[\times\tore).  
 \end{equation*}
 on the pressure, are known in literature as \textit{suitable weak
   solutions} and they are of fundamental importance from the
 theoretical point of view since they are those for which partial
 regularity results holds true, see Scheffer~\cite{Sche1977} and
 Caffarelli-Kohn-Nirenberg~\cite{CKN1982}.

 Due to possible non-uniqueness of solutions in the 3D case, see in
 particular the recent result in~\cite{ABC2021} for the case with
 external forces, it is not ensured that all schemes produce weak
 solutions, with the correct global and local balance. Moreover, from
 the applied point of view the local energy inequality is a sort of
 entropy condition and, even if it is not enough to prove uniqueness,
 it seems a natural request to select physically relevant solutions,
 especially for turbulent or convection dominated problems.  For this
 reason it is natural to ask, in view of obtaining accurate
 simulations of turbulent flows, that the above local energy
 inequality has to be satisfied by solutions constructed by numerical
 methods.  The interplay between suitable weak solutions and numerical
 computations of turbulent flows has been emphasized starting from the
 work of Guermond \textit{et al.}~\cite{Gue2008,GOP2004} and a recent
 overview can be found in the monograph~\cite{Ber2021}. In this paper,
 we continue and extend some previous works in~\cite{Ber2018, BS2012,
   BS2014,BS2016} and especially~\cite{BFS2019} to analyze the
 difficulties arising when dealing with full space-time discretization
 with schemes which are of second order in the time variable.
 The aim of this paper is to extend to the case $\theta=1/2$, which
 corresponds to the Crank-Nicolson method and which could not be
 treated directly with the same proofs as in~\cite{BFS2019}.  In
 particular, the case $\theta=1/2$ requires a coupling between the
 space and time mesh-size, which is nevertheless common to other
 second order models. In fact, beside the Crank-Nicolson
 scheme~\eqref{eq:disintro}, studied in Section~\ref{sec:5}, we will
 also consider in the final Section~\ref{sec:6} other schemes
 involving the Adams-Bashforth or the Linear-Extrapolation for the
 convective term.
 
 To set up the problem we consider as in~\cite{Gue2006} two sequences
 of discrete approximation spaces $\{X_{h}\}_{h}\subset H_{\#}^{1} $
 and $\{M_{h}\}_{h}\subset H_{\#}^{1}$ which satisfy --among other
 properties described in Section~\ref{sec:dis}-- an appropriate
 commutator property, see Definition~\ref{def:dcp}. Then, given a net
 $t_m:=m\,\Delta t$ we consider the following implicit space-time
 discretization of the problem~\eqref{eq:nse}-\eqref{eq:nsid}: Set
 $u_h^0=\pi_h(u_0)$, where $\pi_h$ is the projection over $ X_h$. For
 any $m=1,...,N$ and given $u_h^{m-1}\in X_h$ and
 $p_h^{m-1}\in M_{h}$, find $u_{h}^m\in X_h$ and $p^m\in M_{h}$ such
 that
\begin{equation}
  \label{eq:disintro}
  \tag{CN}  \begin{aligned}
    \left(d_t u^m_h,v_h\right)+\nu(\nabla u_{h}^{m,1/2},\nabla
    v_h)+b_h(u_{h}^{m,1/2}, u_{h}^{m,1/2},v_h)-( p_{h,}^{m},\dive v_h)=0,
    \\
    ( \dive u_{h}^{m},q_h)=0,
  \end{aligned}
\end{equation}
where $d_tu^m:=\frac{u_h^{m}-u_h^{m-1}}{\Delta t}$ is the backward
finite-difference approximation for the time-derivative in the
interval $(t_{m-1},t_m)$ of constant length $\Delta t$;
$\uth:=\frac{1}{2}\big( u^m_h+u^{m-1}_h\big)$ is the average of values
at consecutive time-steps; $b_h(u_{h}^{m,1/2}, u_{h}^{m,1/2},v_h)$ is
a suitable discrete approximation of the non-linear term.  Other
notations, definitions, and properties regarding~\eqref{eq:disintro}
will be given in Sections~\ref{sec:pre}-\ref{sec:dis}.  We refer to
Quarteroni and Valli~\cite{QV1994} and Thom\'ee~\cite{Tho1997} for
general properties of $\theta$-schemes (not only for $\theta=1/2$) for
parabolic equations. Recall that for the fully implicit Crank–Nicolson
scheme, Heywood and Rannacher~\cite{HR1990} proved that it is almost
unconditionally stable and convergent. For a two-step scheme with a
semi-implicit treatment for the nonlinear term, He and
Li~\cite{HL1996} gave the convergence condition:
$\Delta t\, h^{-1/2}\leq C_{0}$. For the
Crank–Nicolson/Adams–Bashforth scheme in which the nonlinear term is
treated explicitly, Marion and Temam~\cite{MT1998} provided the
stability condition $\Delta t \,h^{-2}\leq C_{0}$ and
Tone~\cite{Ton2004} proved the convergence under the condition
$\Delta t \,h^{-2-3/2}\leq C_{0}$, and in all cases
$C_{0}=C_{0}(\nu,\Omega,T,u_{0},f)$. The situation is different in two
space dimensions, cf.~He and Sun~\cite{HS2007}, where more regularity
of the solution can be used, but these results are not applicable to
genuine (turbulent) weak solutions in the three dimensional case.  We
observe that the value $\theta=1/2$ makes the scheme more accurate in
the time variable, but on the other hand introduces some ``natural''
or at least expected limitation on the mesh-sizes. Other schemes will
be also considered, in order to adapt the results also to different
second order schemes, since the proof is rather flexible to handle
several different discretizations of the NSE.

As usual in time-discrete problem (see for instance~\cite{Tem1977b}),
in order to study the convergence to the solutions of the continuous
problem it is useful to consider $v^{\Delta t}_{h}$ which is the
linear interpolation of $\{u^{m}_h\}_{m=1}^{N}$ (over the net
$t_m=m\Delta t$), and $u^{\Delta t}_{h}$ and $p^{\Delta t}_{h}$ which
are the time-step functions such that on the interval
$[t_{m-1}, t_{m})$ are equal to $u^{m,1/2}_h$ and $p^{m}_h$,
respectively, see~\eqref{eq:vm}.

The main result of the paper is the following, we refer to
Section~\ref{sec:pre} for further details on the notations.
\begin{theorem}
  \label{teo:main} Let the finite element spaces $(X_{h},M_{h})$
  satisfy the discrete commutation property, and the technical
  conditions described in Section~\ref{sec:space discretization}.  Let
  $u_{0}\in H^{1}_{\dive}$ and fix $\Delta t>0$ and $h>0$ such that
        \begin{equation}\label{eq:onehalf}
          \begin{aligned}
            &\frac{\Delta t \|u_0\|_2^3}{\nu\,
              h^{1/2}}=o(1),
          \end{aligned}
        \end{equation}
        Let
        $\{(v^{\Delta t}_{h}, u^{\Delta t}_{h}, p^{\Delta
          t}_{h})\}_{\Delta t,h}$ as in~\eqref {eq:vm}. Then, there
        exists
        \begin{equation*}
          (u,p)\in L^\infty(0,T;L_{\dive}^2)\cap L^2(0,T;H^1_{\dive})\times
          L^{4/3}(0,T;L^{2}_{\#}), 
        \end{equation*}
        such that, up to a sub-sequence, as $(\Delta t,h)\to(0,0)$,
        \begin{equation*}
          \begin{aligned}
            &v^{\Delta t}_{h}\rightarrow u\textrm{ strongly in
            }L^{2}((0,T)\times\tore), 
            \\ 
            &u^{\Delta t}_{h}\rightarrow u\textrm{ strongly in
            }L^{2}((0,T)\times\tore), 
            \\ 
            &\nabla u^{\Delta t}_{h}\weakto \nabla u\textrm{ weakly in
            }L^{2}((0,T)\times\tore),
            \\
            &p^{\Delta t}_{h}\weakto p\textrm{ weakly in
            }L^{\frac{4}{3}}((0,T)\times\tore).
          \end{aligned}
        \end{equation*}
        Moreover, the couple $(u,p)$ is a suitable weak solution
        of~\eqref{eq:nse}-\eqref{eq:nsid} in the sense of
        Definition~\ref{def:suit}.
      \end{theorem}
\begin{remark}
  Theorem~\ref{teo:main} holds also in the presence of an external
  force $f$ satisfying suitable bounds. For example,
  $f\in L^{2}(0,T;L^{2}(\tore))$ is enough.
\end{remark}
The proof of Theorem~\ref{teo:main} is given in Section~\ref{sec:5}
and it is based on a compactness argument as we previously developed
in~\cite{BFS2019} and a precise analysis of the quantity
$u^{m+1}_{h}-u^{m}_{h}$, by using the assumptions linking the time and
the spatial mesh size.

\medskip

\textbf{Plan of the paper.}  In Section~\ref{sec:pre} we fix the
notation that we use in the paper and we recall the main definitions
and tools used. In Section~\ref{sec:dis} we introduce and give some
details about the space-time discretization methods. Finally, in
Section~\ref{sec:apriori} we prove the main {\em a priori} estimates
needed to study the convergence and in Section~\ref{sec:5} we prove
Theorem~\ref{teo:main}. In the final Section~\ref{sec:6} we adapt the
proofs to a couple of different second order schemes.
\section{Notations and Preliminaries}
\label{sec:pre}
In this section we fix the notation we will use in the paper; we also 
recall the main definitions concerning weak solutions of
incompressible NSE and a compactness result.
\subsection{Notations}
We introduce the notations typical of space-periodic problems. We will
use the customary Lebesgue spaces $L^p(\tore)$ and Sobolev spaces
$W^{k,p}(\tore)$ and we will denote their norms by $\|\cdot\|_p$ and
$\|\cdot\|_{W^{k,p}}$ We will not distinguish between scalar and
vector valued functions, since it will be clear from the context which
one has to be considered. In the case $p=2$, the $L^{2}(\tore)$ scalar
product is denoted by $\left(\cdot,\cdot\right)$, we use the notation
$H^s(\tore):=W^{s,2}(\tore)$ and we define, for $s>0$, the dual spaces
$H^{-s}(\tore):=(H^s(\tore))'$.  Moreover, we will consider always
sub-spaces of functions with zero mean value and these will be denoted
by
\begin{equation*}
    L_{\#}^{p}:= \left\{w\in L^{p}(\tore):
    \quad\int_{\tore}w\,dx = {0} \right\}\qquad 1\leq p\leq+\infty,
\end{equation*}
and also 
\begin{equation*}
  H_{\#}^{s}:= H^{s}(\tore)\cap L^{2}_{\#}.
\end{equation*}
As usual we consider spaces of divergence free vector fields, defined
as follows
\begin{equation*}
  L_{\dive}^{2}:= \left\{w \in (L^{2}_{\#})^{3}: \ \dive w = 0\right\}
  \qquad\text{and, for $s>0$, }  H_{\dive}^{s} := H^{s}_{\#}\cap
  L_{\dive}^{2}. 
\end{equation*}
Finally, given $X$ a Banach space, $L^p(0,T;X)$ denotes the classical
Bochner spaces of $X$ valued functions, endowed with its natural norm,
denoted by $\|\cdot\|_{L^{p}(X)}$. We denote by $l^p(X)$ the discrete
counterpart for $X$-valued sequences $\{x^m\}$, defined on the net
$\{m\Delta t\}$, and with weighted norm defined by
$\|x\|_{l^{p}(X)}^p:=\Delta t\sum_{m=0}^M\|x^m\|_X^p$.
\subsection{Weak solutions and suitable weak solutions}
We start by recalling the notion of weak solution (as introduced by
Leray and Hopf) and adapted to the space periodic setting.
\begin{definition}
  The vector field $u$ is a Leray-Hopf weak solution
  of~\eqref{eq:nse}-\eqref{eq:nsid} if 
  \begin{equation*}
    u\in L^{\infty}(0,T;L_{\dive}^{2})\cap L^{2}(0,T;H^{1}_{\dive}),
  \end{equation*}
  and if $u$ satisfies the NSE~\eqref{eq:nse}-\eqref{eq:nsid} in the
  weak sense, namely the integral equality
  \begin{equation}
    \label{eq:nsw}
    \int_0^{T}\big[ \left(u,\partial_t\phi\right)-\nu\left(\nabla
      u,\nabla\phi\right)-\left((u\cdot\nabla)\,u,\phi\right)\big]\,dt
    +\left(u_0,\phi(0)\right)=0,   
  \end{equation}
  holds true for all smooth, periodic, and divergence-free functions
  $\phi\in C_c^{\infty}([0,T);C^{\infty}(\tore))$.  Moreover, the
  initial datum is attained in the strong $L^{2}$-sense, that is
  \begin{equation*}
    \lim_{t\to 0^{+}}\|u(t)-u_0\|_{2}=0,
  \end{equation*}
  and the following \textit{global} energy inequality holds
  \begin{equation}
    \label{eq:gei}
    \frac12\|u(t)\|_{2}^2+\nu\int_0^t\|\nabla u(s)\|_{2}^2\,ds\leq\frac12
    \|u_0\|_{2}^2,\quad\textrm{ for all } t\in[0,T].
  \end{equation}
\end{definition}
Suitable weak solutions are a particular subclass of Leray-Hopf weak
solutions and the definition is the following.
\begin{definition}
  \label{def:suit}
  A pair $(u,p)$ is a suitable weak solution to the Navier-Stokes
  equation~\eqref{eq:nse} if $u$ is a Leray-Hopf weak solution,
  $p\in L^{\frac{4}{3}}(0,T;L^{2}_{\#})$, and the local energy
  inequality
  \begin{equation}
  \label{eq:lei}
  \nu\int_0^T\int_{\tore}|\nabla u|^2\phi\,dxdt\leq
  \int_{0}^{T}\int_{\tore}\left[\frac{|u|^2}{2}\left(\partial_t\phi
      +\nu\Delta\phi\right)
    +\left(\frac{|u|^2}{2}+p\right)u\cdot\nabla \phi\right]\,dxdt,
\end{equation}
holds for all $\phi\in C^\infty_0(0,T;C^\infty(\tore))$ such that
$\phi\geq0$, 
\end{definition}
\begin{remark}
  The definition of suitable weak solution is usually stated with
  $p\in L^{\frac53}((0,T)\times\tore)$ while in
  Definition~\ref{def:suit} $p\in L^\frac43(0,T;L^{2}(\tore))$. This
  is not an issue since of course we have a bit less integrability in
  time but we gain a full $L^{2}$-integrability in space. The main
  property of suitable weak solutions is the fact that they satisfy
  the local energy inequality~\eqref{eq:lei} and weakening the
  requests on the pressure does not influence the validity of local
  regularity results, see for instance discussion in
  Vasseur~\cite{Vas2007}.
\end{remark}
\subsection{A compactness lemma}
In this subsection we recall the main compactness lemma which allows
us to prove the strong convergence of the approximations. We remark
that it is a particular case of a more general lemma, whose statement
and proof can be found in~\cite[Lemma 5.1]{Lio1998}. Even if it is a
tool most often used for compressible equations, it is useful here and
we recall the special version taken from~\cite{BFS2019}.
\begin{lemma}
  \label{lem:comp}
  Let $\{f_n\}_{n\in\N}$ and $\{g_n\}_{n\in\N}$ be uniformly bounded
  in $L^{\infty}(0,T;L^{2}(\tore))$ and let be given
  $f,g\in L^{\infty}(0,T;L^{2}(\tore))$ such that
  \begin{equation*}
    \label{eq:lem1}
  \begin{aligned}
    &f_n\weakto f\textrm{ weakly in }L^{2}((0,T)\times\tore),
    \\
    &g_n\weakto g\textrm{ weakly in } L^{2}((0,T)\times\tore).
\end{aligned}
\end{equation*}
Let $p\geq1$ and assume that 
\begin{equation*}
  \label{eq:asslem}
  \begin{aligned}
    &\{\partial_tf_n\}_n\subset L^{p}(0,T;H^{-1}(\tore)),\quad
    &\{g_n\}_{n}
    \subset    L^{2}(0,T;H^{1}(\tore)), 
  \end{aligned}
\end{equation*}
with uniform (with respect to $n\in\N$) bounds on the norms. Then, 
\begin{equation*}
  \label{eq:conc}
  f_n\,g_n\weakto f\,g\textrm{ weakly in }L^{1}((0,T)\times\tore). 
\end{equation*}
\end{lemma}
\section{Setting of the numerical approximation}
\label{sec:dis}
In this section we introduce the space-time discretization of the
initial value problem~\eqref{eq:nse}-\eqref{eq:nsid}. We start by
introducing the space discretization by finite elements.
\subsection{Space discretization}
\label{sec:space discretization}
For the space discretization we strictly follow the setting considered
in~\cite{Gue2006}. Let $\mathcal{T}_{h}$ be a non-degenerate (shape
regular) simplicial subdivision of $\tore$. Let
$\{X_h\}_{h>0}\subset H_{\#}^{1}$ be the discrete space for
approximate velocity and $\{M_h\}_{h>0}\subset L^{2}_{\#}$ be that of
approximate pressure. To avoid further technicalities, we assume as
in~\cite{Gue2006}, that $M_h\subset H^{1}_{\#}$.

We make the following (technical) assumptions on the spaces $X_h$ and $M_h$:
\begin{enumerate}
\item For any $v\in H^{1}_{\#}$ and for any $q\in L^{2}_{\#}$ there exists
  $\{v_h\}_{h}$ and $\{q_{h}\}_{h}$ with $v_h\in X_h$ and $q_h\in M_h$ such that
  \begin{equation}
    \label{eq:apptest}
    \begin{aligned}
      &v_h\to v\quad\textrm{ strongly in }H^{1}_{\#}\quad\textrm{ as }h\to 0,
      \\
      &q_h\to q\quad\textrm{ strongly in }L^2_{\#}\quad\textrm{ as }h\to 0;
    \end{aligned}
  \end{equation}
\item Let $\pi_h : L^2(\tore)\to X_h$ be the $L^2-$projection onto
  $X_h$. Then, there exists $c>0$ independent of $h$ such that,
\begin{equation}
  \label{coer}
  \forall q_h \in M_h \quad \|\pi_h\left(\nabla q_h\right)\|_{2}\geq
  c\|q_h\|_{2}; 
\end{equation}
\item  There is $c$ independent of $h$ such that for all $v\in
  H_{\#}^{1}$ 
  \begin{equation*}
    \begin{aligned}
      &\|v-\pi_h(v)\|_{2}=\inf_{w_h\in X_h}\|v-w_h\|_{2}\leq c\,h\|v\|_{H^1},
      \\
      &\|\pi_h(v)\|_{H^1}\leq c\|v\|_{H^1};
    \end{aligned}
  \end{equation*}
\item There exists $c$ independent of $h$ (inverse inequality) such that
\begin{equation}
  \label{eq:imm}
\forall v_h\in X_h\quad  \|v_h\|_{H^1}\leq c\,h^{-1}\|v_h\|_{2}.
\end{equation} 
\end{enumerate}
Moreover, we assume that $X_h$ and $M_h$ satisfy the following
discrete commutator property.
\begin{definition}
  \label{def:dcp}
  We say that $X_h$ (resp. $M_h$) has the discrete commutator property
  if there exists an operator $P_h\in \mathcal{L}(H^1;X_h)$
  (resp. $Q_h\in \mathcal{L}(L^2;M_h)$) such that for all
  $\phi \in W^{2,\infty}$ (resp. $\phi \in W^{1,\infty}$) and all
  $v_h\in X_h$ (resp. $q_h\in M_h$)
\begin{align}
  & 	\|v_h\phi - P_h(v_h\phi)\|_{H^l}\leq
    c\,h^{1+m-l}\|v_h\|_{H^m}\|\phi\|_{W^{m+1,\infty}},\label{eq:dcp1}
  \\
  & 	\|q_h\phi - Q_h(q_h\phi)\|_{2}\leq
    c\,h\,\|q_h\|_{2}\|\phi\|_{W^{1,\infty}}\label{eq:dcp2}, 
\end{align}
for all $0\leq l \leq m\leq 1$.
\end{definition}
\begin{remark}\label{rem:example}
  Explicit and relevant examples of couples $(X_{h},M_{h})$ of finite
  element spaces satisfying the commutator property are those employed
  in the MINI and Hood-Taylor elements with quasi-uniform mesh,
  see~\cite{Car2002}. 
\end{remark}
We recall from~\cite{Gue2006} that the coercivity
hypothesis~\eqref{coer}, allows us to define the map
$\psi_h : H_{\#}^2 \to M_h$ such that, for all $q\in H_{\#}^2$, the
function $\psi_h(q)$ is the unique solution to the problem:
\begin{equation*}
  \label{eq:preintro}
  \left(\pi_h(\nabla \psi_h(q)),\nabla r_h\right)=\left(\nabla
    q,\nabla r_h\right). 
\end{equation*} 
This map has the following properties: there exists $c$, independent of
$h$, such that for all $q\in H_{\#}^2$, 
\begin{align*}
  & \nonumber\|\nabla (\psi_h(q)-q)\|_{2}\leq c\,h\|q\|_{H^2},
  \\
  &  \|\pi_h\nabla \psi_h(q)\|_{H^1}\leq c\|q\|_{H^2}.\label{eq:psi2}
\end{align*}
Let us introduce the space of discretely divergence-free functions
\begin{equation*}
V_h=\left\{ v_h\in X_h : \left(\dive v_h,q_h\right)=0 \quad \forall
  q_h\in L^2(\Omega)\right\}. 
\end{equation*}

\medskip

The most common variational formulation (for the continuous problems)
of the convective term $nl(u,v):=(u\cdot\nabla)\,v$ is
\begin{equation*}
  b(u,v,w)=\int_{\tore}(u\cdot\nabla)\,v\cdot w\,dx,
\end{equation*}
and the fact that $ b(u,v,v)=0$ for $u\in L^{2}_{\dive}$,
$v\in H^{1}_{\#}$, allows us to deduce, at least formally, the energy
inequality~\eqref{eq:gei}. This cancellation is based on the constraint
$\dive u=0$ and this identity is not valid anymore in the case of
discretely divergence-free functions in $V_{h}$. To have the basic
energy estimate we need to modify the non-linear term since
$V_h\not\subseteq H^{1}_{\dive}$ and the choice of the weak
formulation becomes particularly relevant in the discrete case since
it leads to schemes with very different numerical properties.

To formulate the various schemes we will consider, which corresponds to
the Cases~1-2-3, we define the discrete tri-linear operator
$b_{h}(\cdot,\cdot,\cdot)$ in different (but standard) ways. This
permits a sort of unified treatment: for instance in all the three
case considered below it holds \textit{at least} that $nl_h(u,v)$
--which is the discrete counterpart of $nl(u,v)$-- satisfies the
following estimate
\begin{equation*}
  \|nl_h(u,v)\|_{H^{-1}}\leq \|u\|_{3}\,\|v\|_{H^1}\quad
  \forall\,u,v\in H^{1}_{\#}. 
\end{equation*}
\smallskip

We present now with details the various different discrete
formulations we will use.

\medskip

\noindent\textbf{Case~1}: 
We use the most common option, that is of a ``symmetrized'' operator 
\begin{equation} 
  \label{eq:nonlinear} 
  nl_h(u,v):=(u\cdot\nabla)\, v+ \frac{1}{2}v\dive u,
\end{equation}
for the convective term, which leads to the tri-linear form
\begin{equation}
  \label{eq:trilinear} 
  b_h(u,v,w):=\langle   nl_h(u,v), w\rangle_{H^{-1}\times    H_{\#}^1},
\end{equation}
such that 
\begin{equation*}
  \label{eq:nonlinear2}
  b_h(u,v,v)=0\qquad \forall\,u,\,v\in H^{1}_{\dive}+V_h.
\end{equation*}
Moreover, this tri-linear operator can be also estimated as follows
\begin{equation}
\label{eq:estimate-nonlinear1}
\begin{aligned}
    |b_h(u,v,w)|&\leq \|u\|_{6}\|\nabla
    v\|_{2}\|w\|_{3}+\frac{1}{2}\|v\|_{6}\|\dive u\|_{2}\|w\|_{3} 
    \\
&\leq    C\|\nabla u\|_{2}\|\nabla v\|_{2}\|w\|_{2}^{1/2}\|\nabla
w\|_{2}^{1/2} ,
  \end{aligned}
\end{equation}
by means of the Sobolev embedding $H^{1}(\tore)\subset L^{6}(\tore)$
and of the convex interpolation inequality.

\bigskip

\noindent\textbf{Case~2}:
Alternatively, we can consider  the ``rotational form without pressure,'' as
in Layton \textit{et al.}~\cite{LMNOR2009}, which corresponds to the
formulation
\begin{equation} 
  \label{eq:nonlinear-rotational-pressure} 
  nl_h(u,v):=(\nabla\times u)\times  v,
\end{equation}
and which leads to the tri-linear form 
$$
b_h(u,v,w):=\langle nl_h(u,v), w\rangle_{H^{-1}\times H_{\#}^1},
$$
such that 
\begin{equation*}
  b_h(u,v,v)=0\qquad \forall\,u,\,v\in H^{1}_{\dive}+V_h. 
\end{equation*}
Moreover, this term can be estimated as follows
\begin{equation}
  \label{eq:estimate-nonlinear1bis}
  \begin{aligned}
    |b_h(u,v,w)|&\leq \|\nabla\times u\|_{2}\|v\|_{6}\|w\|_{3}
    \\
    &\leq    C\|\nabla u\|_{2}\|\nabla v\|_{2}\|w\|_{2}^{1/2}\|\nabla
    w\|_{2}^{1/2} ,
  \end{aligned}
\end{equation}
by means again of the Sobolev embedding and of the convex
interpolation inequality. In this case one is \textit{hiding} the
Bernoulli pressure $\frac{1}{2}|v|^{2}$ into the kinematic
pressure. It is well documented that the scheme is easier to be
handled but the {\em under resolution} of the pressure has some
effects on the accuracy, see Horiuti~\cite{Hor1987},
Zang~\cite{Zang1991}, and the discussion in~\cite{LMNOR2009}.

\medskip

In order to overcome the numerical problems arising when using the
operator from Case~2, other computationally more expensive methods are
considered, as the one below

\bigskip

\noindent\textbf{Case~3:}
We consider the rotational form with approximation of the
Bernoulli pressure, as studied already in Guermond~\cite{Gue2006}.
\begin{equation*} 
  nl_h(u,v):=(\nabla\times u)\times  v+ \frac{1}{2}\nabla
  \big(\mathcal{K}_{h}(v\cdot u)\big), 
\end{equation*}
where $\mathcal{K}_{h}$ is the $L^2\to M_{h}$ projection operator,
which is stable, linear, and is defined as
$(\mathcal{K}_{h}u,v_{h})=(u,v_{h})$, for all $u\in L^{2}$ and
$v_{h}\in M_{h}$. In this way the tri-linear term is such that
\begin{equation*}
b_h(u,v,w):=\langle nl_h(u,v), w\rangle_{H^{-1}\times H_{\#}^1},
\end{equation*}
such that 
\begin{equation*}
  b_h(u,v,v)=0\qquad \forall\,u,\,v\in H^{1}_{\dive}+V_h.
\end{equation*}
A first estimate, which is proved also in~\cite{Gue2006}, is the
following one:
\begin{equation*}
 | b_{h}(v,v,w)|\leq c \|v\|_{H^{1}}\|v\|_{{3}}\,\|w\|_{H^{1}}.
\end{equation*}
Here, to better estimate the effect of the projection of the Bernoulli
pressure, we use some improved properties of the $L^{2}$-projection
operator $\mathcal{K}_{h}$, which are valid in the case of
quasi-uniform meshes. In fact, for special meshes one can show also
the $W^{1,p}$ stability. The improved stability for the
$L^{2}$-projection has a long history, see Douglas, Jr., Dupont, and
Wahlbin~\cite{DDW1974}, Bramble and Xu~\cite{BX191}, and the review in
the recent work of Diening, Storn, and Tscherpel~\cite{DST2021}. A
different approach in Hilbert fractional spaces is used in
\cite{Gue2007}, but we do not known whether this applies to the
estimates we are willing to use.  Anyway, we couldn't find the
detailed proof of the required stability, which can be obtained by
using the $L^{p}$-stability of the operator $\mathcal{K}_{h}$, the
inverse inequality valid for the meshes we consider, and the
$W^{1,p}$-stability and approximation of the Scott-Zhang projection
operator ${\Pi}_{h}$ (see~\cite{BS2008}) valid for
$f\in W^{1,p}(\Omega)$, for $p\in[1,\infty[$. Just to sketch the
argument, it is enough to use the following inequalities
\begin{equation*}
  \begin{aligned}
    \| \mathcal{K}_{h}f \|_{W^{1,p}}&\le c
    \| \mathcal{K}_{h}(f-\Pi_{h}f) \|_{W^{1,p}}
    + \| \Pi_{h} f \|_{W^{1,p}} 
    \\
    &\le c h^{-1} \| \mathcal{K}_{h}(f-\Pi_{h}f) \|_p
    + \| \Pi_{h} f \|_{W^{1,p}}
    \\
    &\le c h^{-1} \| f-\Pi_{h} f
    \|_{p} + \| \Pi_{h}f \|_{W^{1,p}}
    \\
    & \le c \|f\|_{W^{1,p}}.
  \end{aligned}
\end{equation*}
With this stability result, in the ``case 3'' the nonlinear term can
be estimated as follows
\begin{equation}
\label{eq:estimate-nonlinear1tris}
\begin{aligned}
    |b_h(u,v,w)|&\leq \|\nabla u\|_{2}\|
    v\|_{6}\|w\|_{3}+\|\nabla\mathcal{K}_{h}(u\cdot v)\|_{3/2}\|w\|_{3}
    \\
&\leq    \big(C\|\nabla u\|_{2}\|\nabla
v\|_{2}+\|\mathcal{K}_{h}(u\cdot v)\|_{W^{1,3/2}}\big)\|w\|_{3} 
\\
&\leq \big(C\|\nabla u\|_{2}\|\nabla v\|_{2}+\|u\|_{3}\| v\|_{3}+\|\nabla
u\|_{2}\| v\|_{6}+\|u\|_{6}\| \nabla v\|_{2}\big)\|w\|_{3} 
\\
& \leq C\|\nabla u\|_{2}\|\nabla v\|_{2}\|w\|_{2}^{1/2}\|\nabla w\|_{2}^{1/2},
\end{aligned}
\end{equation}
by means of the Sobolev embedding $H^{1}(\tore)\subset L^{6}(\tore)$
and of the convex interpolation inequality, exactly as in the first
two cases.
\subsection{Time discretization}
We now pass to the description of the time discretization. For the
time variable $t$ we define the mesh as follows: Given $N\in\N$ the
time-step $0<\Delta t\leq T$ is defined as $\Delta
t:=T/N$. Accordingly, we define the corresponding net
$\{t_m\}_{m=1}^{N}$ by
\begin{equation*}
  t_0:=0\,\qquad
  t_m:=m\,\Delta t,\qquad m=1,\dots,N.
\end{equation*}
We consider the (Crank-Nicolson) method~\eqref{eq:disintro}
(cf.~\cite[\S~5.6.2]{QV1994}). With a slight abuse of notation we
consider $\Delta t=T/N$ and $h$, instead of $(N,h)$, as the indexes of
the sequences for which we prove the convergence. Then, the
convergence will be proved in the limit as $(\Delta t,h)\to(0,0)$. We
stress that this does not affect the proofs since all the convergences
are proved up to sub-sequences.\par
Once~\eqref{eq:disintro} is solved, we consider a continuous version
useful to study the convergence. To this end we associate to the
triple $(\uth,u_{h}^{m},p^m_{h})$ the functions
\begin{equation*}
(v^{\Delta t}_{h},u^{\Delta t}_{h}, p^{\Delta
  t}_{h}):[0,T]\times\tore\rightarrow \R^3\times\R^3\times \R, 
\end{equation*}
defined as follows:
 \begin{equation}
 \label{eq:vm}
 \begin{aligned}
   &v^{\Delta t}_{h}(t):=\left\{
     \begin{aligned}
       &u_{h}^{m-1}+\frac{t-t_{m-1}}{\Delta t}(u_{h}^{m}-u_{h}^{m-1})\ & \text{for
       } t\in[t_{m-1},t_m),
       \\
       &u_{h}^{N}\qquad &\text{for
       } t=t_N,
     \end{aligned}
   \right.
   \\
   &u^{\Delta t}_{h}(t):=\left\{
     \begin{aligned}
       &\uth & \text{for
       } t\in[t_{m-1},t_m),
       \\
       &u_{h}^{N,1/2}\qquad &\text{for
       } t=t_N,
     \end{aligned}
   \right.
   \\
      &p^{\Delta t}_{h}(t):=\left\{
     \begin{aligned}
       &p_h^{m} & \text{for
       } t\in[t_{m-1},t_m),
       \\
       &p_{h}^{N}\qquad &\text{for
       } t=t_N.
     \end{aligned}
   \right.
\end{aligned}
\end{equation}
Then, the discrete equations~\eqref{eq:disintro} can be rephrased as
the following time-continuous system:
\begin{equation}
  \label{eq:thet1ac}
  \begin{aligned}
    \left(\partial_t v^{\Delta t}_{h},w_h\right)+b_h\left(u^{\Delta
        t}_{h},u^{\Delta t}_{h},w_h\right)+\nu\left(\nabla u^{\Delta
        t}_{h},\nabla w_h\right)-\left(p^{\Delta t}_{h},\dive
      w_h\right)&=0,
    \\
    \left(\dive u^{\Delta t}_{h}, q_h\right)&=0,
  \end{aligned}
\end{equation}
for all $w_h\in L^{s}(0,T;X_h)$ (with $s\geq4)$ and for all
$q_h\in L^{2}(0,T;M_h)$. We notice that the divergence-free condition
comes from the fact that $u^{m}_{h}$ is such that
\begin{equation*}
  \left(\dive u^m_h, q_{h}\right)=0\qquad\text{ for }m=1,...,N,\
  \forall\,q_{h}\in M_{h}. 
\end{equation*}
\section{A priori estimates}
\label{sec:apriori}
In this section we prove the {\em a priori} estimates that we need to
study the convergence of solutions of~\eqref{eq:thet1ac} to suitable
weak solutions of~\eqref{eq:nse}-\eqref{eq:nsid}. We start with the
following discrete energy equality.
\begin{lemma}
  \label{lem:discene}
  Let $N\in\N$ and $m=1,..,N$. Then, for solutions
  to~\eqref{eq:disintro} the following (global) discrete energy-type
  equality holds true:
  \begin{align}
    \label{eq:abis}
    & \frac{1}{2}(\|u_h^m\|_{2}^2-\|u_h^{m-1}\|_{2}^2)+\nu\Delta t \|\nabla
    u_{h}^{m,\frac{1}{2}}\|_{2}^2=0.
  \end{align}
  Moreover, if $u_{0}\in H^{1}_{\dive}$ there exists $C>0$ depending on
  $\|u_{0}\|_{H^{1}}$ such that
  \begin{equation}
    \label{eq:atris}
    \sum_{m=1}^{N}    \|u_h^m-u_h^{m-1}\|_{2}^2\leq C\left(\Delta
      t+\frac{1}{h^{1/2}}\right).
  \end{equation}
\end{lemma}
\begin{proof}
  We start by proving the (global) discrete energy equality. For any
  $m=1,...,N$ take  $w_h=\chi_{[t_{m-1},t_m)}u_{h}^{m,1/2}\in
  L^{\infty}(0,T;X_h)$ as test function in~\eqref{eq:thet1ac}. Then,
  it follows 
  \begin{equation*}
    \left(\frac{u_h^m-u_h^{m-1}}{\Delta
        t},u_{h}^{m,1/2}\right)+\nu\|\nabla u_{h}^{m,1/2}\|_{2}^2=0, 
  \end{equation*}
  which holds true since $u_{h}^{m,1/2}\in X_{h}$ and
  $p^{m}_{h}\in M_{h}$, we have that
  \begin{equation*}
    b_h(u_{h}^{m,1/2}, u_{h}^{m,1/2}, u_{h}^{m,1/2})=0\qquad
    \text{and}\qquad (p_{h}^{m},\dive    u_{h}^{m,1/2})=0.  
  \end{equation*}
  The term involving the discretization of the time-derivative reads
  as follows:
  \begin{align*}
    (u_h^m-u_h^{m-1},u_{h}^{m,1/2})=
    \frac{1}{2}(u_h^m-u_h^{m-1},u_h^m+u_h^{m-1}) 
       =\frac{1}{2}(\|u_h^m\|_2^2-\|u_h^{m-1}\|_2^2).
  \end{align*}
  Then, multiplying by $\Delta t>0$, Eq.~\eqref{eq:abis} holds
  true. In addition, summing over $m$ we also get
  \begin{equation*}
    \frac{1}{2}\|u_h^N\|_2^2
    + \nu\Delta    t\,\sum_{m=0}^N\|\nabla   u_{h}^{m,1/2}\|_2^2=
    \frac{1}{2}\|u_h^0\|_2^2, 
  \end{equation*}
  which proves the $l^\infty(L^2_{\#})\cap l^2(H^1_\#)$ uniform bound
  for the sequence $\{u^{m}_{h}\}$.

  To prove~\eqref{eq:atris} take
  $w_h=\chi_{[t_{m-1},t_m)}(u_{h}^{m}-u_{h}^{m-1})\in
  L^{\infty}(0,T;X_h)$ in~\eqref{eq:thet1ac}. Then, after
  multiplication by $\Delta t$ we get
  \begin{equation*}
    \begin{aligned}
      \|u_{h}^{m}-u_{h}^{m-1}\|_{2}^{2}&+\frac{\nu \Delta t}{2}
      (\|\nabla u_{h}^{m}\|^{2}-\|\nabla u_{h}^{m-1}\|^{2}) 
      \\
      &\leq \Delta t|b_{h}(
      u_{h}^{m,\frac{1}{2}},u_{h}^{m,\frac{1}{2}},u_{h}^{m}-u_{h}^{m-1})| 
      \\
      &\leq C \Delta
      t\|\nabla
      u_{h}^{m,\frac{1}{2}}\|_{2}^{2}\|u_{h}^{m}-u_{h}^{m-1}\|_{2}^{1/2}
      \|\nabla(u_{h}^{m}-u_{h}^{m-1})\|_{2}^{1/2},  
    \end{aligned}
  \end{equation*}
  where in the last line~\eqref{eq:estimate-nonlinear1} has been
  used. By using the inverse inequality~\eqref{eq:imm} and summing
  over $m=1,\dots,N$ we get
  \begin{equation*}
    \begin{aligned}
      \frac{\nu \Delta t}{2}\|\nabla u_{h}^{N}\|^{2}_{2}+\sum_{m=1}^{N}
      \|u_{h}^{m}-u_{h}^{m-1}\|_{2}^{2}&\leq \frac{\nu \Delta
        t}{2}\|\nabla u_{h}^{0}\|^{2}_{2} +C\frac{\Delta
        t}{h^{1/2}}\sum_{m=1}^{N}\|\nabla
      u_{h}^{m,\frac{1}{2}}\|_{2}^{2}\|u_{h}^{m}-u_{h}^{m-1}\|_{2}
      \\
      &\leq \frac{\Delta t}{2}\|\nabla u_{h}^{0}\|^{2}_{2}+\frac{\Delta
        t}{h^{1/2}}\sum_{m=1}^{N}\|\nabla
      u_{h}^{m,\frac{1}{2}}\|_2^{2}(\|u_{h}^{m}\|_{2}+\|u_{h}^{m-1}\|_{2})
      \\
      &\leq \frac{\Delta t}{2}\|\nabla
      u_{h}^{0}\|^{2}_{2}+2C\|u_{h}^{0}\|_{2} \frac{\Delta
        t}{h^{1/2}}\sum_{m=1}^{N}\|\nabla
      u_{h}^{m,\frac{1}{2}}\|_2^{2}
      \\
      & \leq \frac{\Delta t}{2}\|\nabla u_{h}^{0}\|^{2}_{2}+
      \frac{2C\|u_{h}^{0}\|_{2}^{3}}{\nu\, h^{1/2}},
    \end{aligned}
  \end{equation*}
  where we used the $l^\infty(L^2_{\#})\cap l^2(H^1_\#)$ bounds coming
  from the energy equality, and then proving the thesis.
\end{proof}
\begin{remark}
  At first glance the inequality~\eqref{eq:atris} seems useless, being
  badly depending on $h$. Recall that the convergence to zero of
  $\sum_{m=1}^{N} \|u_{h}^{m}-u_{h}^{m-1}\|_{2}^{2}$ is a required
  step to identify the limits of $v^{\Delta t}_{h}$ and of
  $u^{\Delta t}_{h}$. Nevertheless, the key step in the next section
  will be that of combining this inequality with a standard
  restriction on the ratio between time and space mesh-size, to
  enforce the equality of the two limiting functions.
\end{remark}
The next lemma concerns the regularity of the pressure. We follow the
argument in~\cite{Gue2006} and we notice that we are essentially
solving the standard discrete Poisson problem associated to the
pressure. It is for this result that the space-periodic setting is
needed.
\begin{lemma}
  \label{lem:pres}
  There exists a constant $c>0$, independent of $\Delta t$ and of $h$,
  such that
  \begin{equation*}
  \begin{aligned}
    \|\pth\|_{2} &\leq
    c\left(\|u_h^{m,1/2}\|_{H^1}+\|u_h^{m,1/2}\|_{3}\,
      \|u_h^{m,1/2}\|_{H^1}\right)\qquad \text{for }m=1,\dots,N.
  \end{aligned}
\end{equation*}
\end{lemma}
\begin{proof}
The proof is exactly the same as in~\cite[Lemma~4.3]{BFS2019}.
\end{proof}
We are now in position to prove the main a priori estimates on the
approximate solutions of~\eqref{eq:thet1ac}.
    \begin{proposition}
      \label{prop:1}
      Let $u_0\in L^{2}_{\dive}$ and assume that~\eqref{eq:onehalf}
      holds. Then, there exists a constant $c>0$, independent of
      $\Delta t$ and of $h$, such that
      \begin{itemize}
      \item[a)] $\|v^{\Delta t}_{h}\|_{L^\infty(L^2)}\leq c$,
      \item[b)]
        $ \|u^{\Delta t}_{h}\|_{L^\infty(L^2)\cap L^2(H^1)}\leq c$,
      \item[c)] $\|p^{\Delta t}_{h}\|_{L^{4/3}(L^{2})}\leq c$,
      \item[d)] $\|\partial_t v^{\Delta t}_{h}\|_{L^{4/3}(H^{-1})}\leq c$,
      \end{itemize}
      Moreover, we also have the following estimate
      \begin{equation}
        \label{eq:newest1}
        \int_0^T\|u^{\Delta t}_{h}-v^{\Delta t}_{h}\|_2^2\,dt\leq
        \frac{\Delta  t}{12}\sum_{m=1}^N\|u_h^m-u_h^{m-1}\|_{2}^2.  
      \end{equation}
    \end{proposition}
    \begin{proof}
      The bound in
      $ L^{\infty}(0,T;L^{2}_{\#})\cap L^{2}(0,T;H^{1}_{\#})$ for
      $v^{\Delta t}_{h}$ follows from~\eqref{eq:vm} and
      Lemma~\ref{lem:discene}, as well as the bounds on
      $u^{\Delta t}_{h}$ in b). The bound on the pressure
      $p^{\Delta t}_{h}$ follows again from~\eqref{eq:vm} and
      Lemma~\ref{lem:pres}. Finally, the bound on the time derivative
      of $v^{\Delta t}_{h}$ follows by~\eqref{eq:thet1ac} and a
      standard comparison argument. Concerning~\eqref{eq:newest1}, by
      using the definitions in~\eqref{eq:vm} we get for
      $t\in[t_{m-1}, t_m)$
      \begin{align*}
        u^{\Delta t}_{h}-v^{\Delta t}_{h}
        &=\frac{1}{2}\,u_h^m+(1-\frac{1}{2})\,u_h^{m-1}-u_{h}^{m-1}
          -\frac{t-t_{m-1}}{\Delta 
          t}(u_{h}^{m}-u_{h}^{m-1})
        \\
        & =\left(\frac{1}{2}-\frac{t-t_{m-1}}{\Delta
          t}\right)\left(u_h^{m}-u_{h}^{m-1}\right).
      \end{align*}
      Then, evaluating the integrals, we have
      \begin{align*}
        \int_0^T\|u^{\Delta t}_{h}-v^{\Delta t}_{h}\|_{2}^2\, dt
        &=\sum_{m=1}^N\|u_h^{m}-u_{h}^{m-1}\|_{2}^2
          \int_{t_{m-1}}^{t_m}\left(\frac{1}{2}-\frac{t-t_{m-1}}{\Delta
          t}\right)^2 dt
        \\
        & \leq\frac{\Delta t}{12}\sum_{m=1}^N\|u_{h}^m-u_{h}^{m-1}\|_{2}^2.
      \end{align*}
    \end{proof}
\section{Proof of the main theorem}
\label{sec:5}
In this section we prove Theorem~\ref{teo:main}. We split the proof in
two main steps: The first one concerns showing that discrete solutions
converge to a Leray-Hopf weak solution, while the second consists in
proving that the constructed solutions are in fact suitable.

\begin{proof}[Proof of Theorem~\ref{teo:main}] We first prove the
  convergence of the numerical sequence to a Leray-Hopf weak solution,
  mainly proving the correct balance of the global
  energy~\eqref{eq:gei}; then, we prove that the weak solution
  constructed is suitable, namely that it satisfies the local energy
  inequality~\eqref{eq:lei}.

\medskip

\paragraph{\textbf{Step 1: Convergence towards a Leray-Hopf weak
    solution}}
We start by observing that from a simple density argument, the test
functions considered in~\eqref{eq:nsw} can be chosen in the space
$L^{s}(0,T;H^{1}_{\dive})\cap C^1(0,T;L^{2}_{\dive})$, with
$s\geq4$. In particular, by using~\eqref{eq:apptest} for any
$w\in L^{s}(0,T;H^{1}_{\dive})\cap C^1(0,T;L^{2}_{\dive})$ such that
$w(T,x)=0$ we can find a sequence
$\{w_h\}_{h}\subset L^{s}(0,T;H^{1}_{\#})\cap C(0,T;L^{2}_{\#})$ such
that
  \begin{equation}
    \label{eq:3}
    \begin{aligned}
      &w_h\to w\textrm{ strongly in
      }L^{s}(0,T;H^{1}_{\#})\quad\textrm{ as }h\to 0,
      \\
      &w_h(0)\to w(0)\textrm{ strongly in }L^{2}_{\#}\quad\textrm{ as
      }h\to 0,
      \\
      &\partial_tw_h\weakto\partial_tw\textrm{ weakly in
      }L^{2}(0,T;L^{2}_{\#})\quad\textrm{ as }h\to 0.
    \end{aligned}
  \end{equation}
  Let $\{(v^{\Delta t}_{h}, v^{\Delta t}_{h}, p^{\Delta t}_{h})\}_{(\Delta t,h)}$, defined
  as in~\eqref{eq:vm}, be a family of solutions of~\eqref{eq:thet1ac}. By
  Proposition~\ref{prop:1}-a) we have that
  \begin{align*}
    & \left\{v^{\Delta t}_{h}\right\}_{(\Delta t,h)}\subset
      L^\infty(0,T;L^{2}_{\#}),\quad\text{ 
      with uniform bounds on the norms}.
  \end{align*}
  Then, by standard compactness arguments there exists $v\in
  L^{\infty}(0,T;L^{2}_{\#})$,  such that (up to a sub-sequence) 
  \begin{equation}
    \label{eq:convv}
    \begin{aligned}
      &v^{\Delta t}_{h}\weakto v\textrm{ weakly in
      }L^{2}(0,T;L^{2}_{\#})\quad \textrm{ as }(\Delta t,h)\to(0,0).
    \end{aligned}
  \end{equation}
  Again by using Proposition~\ref{prop:1} b), there exists $u\in
  L^{\infty}(0,T;L^{2}_{\#})$ such that (up to a sub-sequence)
  \begin{equation}
    \label{eq:convu}
    \begin{aligned}
      &u^{\Delta t}_{h}\overset{*}{\weakto} u\textrm{ weakly* in
      }&L^{\infty}(0,T;L^{2}_{\#})\quad&&\textrm{ as 
      }(\Delta t,h)\to(0,0), 
      \\
      &u^{\Delta t}_{h} \weakto u\textrm{ weakly in }&
      L^2(0,T;H^{1}_{\#})\quad &&\textrm{ as 
      }(\Delta t,h)\to(0,0).
    \end{aligned}
  \end{equation} 
  Moreover, by using~\eqref{eq:apptest}, for any
  $q\in L^2(0,T;L^{2}_{\#})$ we can find a sequence
  $\{q_h\}_{h}\subset L^{2}(0,T;L^{2}_{\#})$ such that
  $q_h\in L^{2}(0,T;M_h)$ and
  \begin{equation*}
    q_h\to q\textrm{ strongly in }L^{2}(0,T;L^{2}_{\#})\quad \textrm{
      as }h\to 0. 
  \end{equation*}
  Then, by using~\eqref{eq:convu} and~\eqref{eq:thet1ac} we have that 
  \begin{equation*}
    0=\int_0^T\left(\dive u^{\Delta t}_{h},q_h\right)\,dt\to \int_0^T\left(\dive
      u,q\right)\,dt\quad \textrm{ as }(\Delta t,h)\to(0,0),
  \end{equation*}
  hence $u$ is divergence-free, since it belongs to
  $H^{1}_{\dive}$. Let us consider~\eqref{eq:newest1}, then
  \begin{equation}
    \label{eq:2}
    \int_{0}^{T}\|v^{\Delta t}_{h}-u^{\Delta t}_{h}\|_{2}^2\,dt\leq\frac{\Delta t}{12}
    \sum_{m=1}^N\|u_h^m-u_h^{m-1}\|_{2}^2 \leq C\,\bigg((\Delta
    t)^{2}+\frac{\Delta t}{h^{1/2}}\bigg),  
  \end{equation}
  where in the last inequality we used again Proposition~\ref{prop:1}
  and the estimate~\eqref{eq:atris}. Then, the integral
  $\int_{0}^{T}\|v^{\Delta t}_{h}-u^{\Delta t}_{h}\|_{2}^2\,dt$
  vanishes as $\Delta t\to0$ if $\Delta t=o(h^{1/2})$, that is
  if~\eqref{eq:onehalf} is satisfied. Then, by using~\eqref{eq:convv}
  and~\eqref{eq:convu} it easily follows that $v=u$.

The rest of the proof follows as in~\cite{BFS2019}, hence we just
sketch the proof, referring to that reference for full details.

By Lemma~\ref{lem:comp} and the fact that $u=v$ we get that
$ u^{\Delta t}_{h}\,v^{\Delta t}_{h}\weakto|u|^{2}$
$\textrm{ weakly in }L^{1}((0,T)\times\tore)$
$\textrm{ as }(\Delta t,h)\to(0,0)$.   In particular, by
using~\eqref{eq:2} we have that
  \begin{equation}
    \label{eq:convuv}
    \begin{aligned}
      &v^{\Delta t}_{h},\, u^{\Delta t}_{h}\to u\textrm{ strongly in
      }L^{2}(0,T;L^{2}_{\#})\quad\textrm{ as 
      }(\Delta t,h)\to(0,0).
    \end{aligned}
  \end{equation}
  Concerning the pressure term the uniform bound in
  Proposition~\ref{prop:1} d) ensures the existence of
  $p\in L^{\frac{4}{3}}(0,T;L^{2}_{\#})$ such that (up to a
  sub-sequence)
  \begin{equation}
    \label{eq:convp}
    p^{\Delta t}_{h}\weakto p\textrm{ weakly in
    }L^{\frac{4}{3}}(0,T;L^{2}_{\#})\quad \textrm{ as 
    }(\Delta t,h)\to(0,0).  
  \end{equation}
  Then, by using~\eqref{eq:3} and~\eqref{eq:convv} we have that
  \begin{equation*}
    \begin{aligned}
      \lim_{(\Delta t,h)\to(0,0)}\int_{0}^{T}(\partial_tv^{\Delta
        t}_{h},w_h)\,dt&
      =-\int_{0}^{T}(u,\partial_tw)\,dt-(u_0, w(0)),
    \end{aligned}
  \end{equation*}
  Next, by using~\eqref{eq:convu},~\eqref{eq:3}, and~\eqref{eq:convp} we
  also get
  \begin{equation*}
    \begin{aligned}
      \lim_{(\Delta t,h)\to(0,0)}\int_{0}^{T}(\nabla u^{\Delta
        t}_{h},\nabla w_{h})\,dt=\int_{0}^{T}(\nabla u, \nabla w)\,dt.
      \\
      \int_{0}^{T}(p^{\Delta t}_{h},\dive w_h)\,dt\to 0\quad \textrm{
        as }(\Delta t,h)\to(0,0).
  \end{aligned}
  \end{equation*}
  Concerning the non-linear term, let $s\geq4$,
  with a standard compactness argument
  \begin{equation*}
    nl_h\left(u^{\Delta t}_{h},u^{\Delta t}_{h}\right) \rightharpoonup
    u\cdot\nabla u,     \mbox{ in } L^{s'}(0,T;H^{-1})\quad \textrm{
      as }(\Delta t,h)\to(0,0).  
  \end{equation*}
  Then, by using also~\eqref{eq:3} it follows that
  \begin{equation}
    \label{eq:convergence-convective}
    \int _0^T b_h(u^{\Delta t}_{h},u^{\Delta t}_{h},w_h)\, dt \to \int _0^T
    \big((u\cdot\nabla)\,u,w\big)\, dt\quad\textrm{ as }(\Delta
    t,h)\to(0,0).  
  \end{equation}
  Finally, the energy inequality follows by Lemma~\ref{lem:discene},
  by using the lower semi-continuity of the $L^{2}$-norm with respect
  to the weak convergence, since the estimate~\eqref{eq:abis} can be
  rewritten as
  \begin{equation*}
    \frac{1}{2}\|v^{\Delta t}_{h}(T)\|^{2}_{2}+\nu\int_{0}^{T}\|\nabla
    u^{\Delta t}_{h}(t)\|^{2}_{2}\,dt\leq\frac{1}{2}\|u_{0}\|^{2}_{2}.
  \end{equation*}

  \medskip
  
  The treatment of the Case~2 and Case~3 can be done with minor
  changes, concerning the tri-linear term, just using the estimate
  already proved in the previous section. The other terms are
  unchanged and the energy estimate remains the same.
  \begin{remark}
    The results for Case~2 and Case~3 can be easily adapted also to
    cover the $\theta$-scheme for $\theta>1/2$, hence completing the
    results in~\cite{BFS2019} which were focusing only on the
    treatment of Case~1.
\end{remark}
Note that, the tri-linear term based on the rotational formulation from
Case~2 and Case~3
$(\nabla\times u^{\Delta t}_{h})\times u^{\Delta t}_{h}$ converges
exactly as in the previous step, since
  \begin{equation*}
      (\nabla\times u^{\Delta t}_{h})\times u^{\Delta t}_{h}
    \rightharpoonup (\nabla\times  u)\times u,
    \mbox{ in } L^{s'}(0,T;H^{-1})\quad \textrm{ as }(\Delta t,h)\to(0,0),
  \end{equation*}
  which implies~\eqref{eq:convergence-convective}, ending the proof in
  the Case~2.

  In the Case~3 the term which needs some care is the projected
  Bernoulli pressure. In this case note that, for
  $\frac{1}{s^{*}}+\frac{1}{2}=\frac{1}{s'}$
\begin{equation*}
  \begin{aligned}
    \int_{0}^{T}\|\mathcal{K}_{h}(|u^{\Delta
      t}_{h}|^{2})-\mathcal{K}_{h}(|u|^{2})\|^{s'}_{2}\,dt &\leq 
    \int_{0}^{T}\||u^{\Delta t}_{h}|^{2}-|u|^{2}\|^{s'}_{2}\,dt
  \\
    &
    \leq  \int_{0}^{T}\||u^{\Delta t}_{h}-u|\, |u^{\Delta
      t}_{h}+u|\|^{s'}_{2}\,dt
    \\
    &\leq\|u^{\Delta t}_{h}-u\|_{L^{s^{*}}(L^{3})}(\|u^{\Delta
      t}_{h}\|_{L^{2}(L^{6})}+\|u\|_{L^{2}(L^{6})}).
  \end{aligned}
\end{equation*}
This shows that
\begin{equation*}
  \mathcal{K}_{h}(|u^{\Delta t}_{h}|^{2})\to \mathcal{K}_{h}(|u|^{2})
  \qquad\text{in }L^{s'}(0,T;L^{2}(\tore)).
\end{equation*}
Moreover, since $\mathcal{K}_{h}(|w|^{2})\to |w|^{2}$ in
$L^{2}(\tore)$ for a.e. $t\in(0,T)$ and
$\|\mathcal{K}_{h}(|w|^{2})\|_{2}\leq \|w\|^{2}_{2}$, by Lebesgue
dominated convergence we have $\mathcal{K}_{h}(|w|^{2})\to |w|^{2}$ in
$L^{s'}(0,T;L^{2}(\tore))$, finally showing that
\begin{equation*}
  \int_{0}^{T}(\mathcal{K}_{h}(|u^{\Delta t}_{h}|^{2}),\dive v)\,dt\to
  \int_{0}^{T}(|u|^{2},\dive v)\,dt.
\end{equation*}
This proves, after integration by parts, that 
\begin{equation*}
  \int_{0}^{T}b_{h}(u^{\Delta t}_{h},
  u^{\Delta t}_{h},w)\,dt \to\int_{0}^{T}(u\cdot\nabla) \,u,w)\,dt.
\end{equation*}
Note also that since in all cases it holds that
\begin{equation}
  \label{eq:energy-inequality-approximated}
  \frac{1}{2}\|v^{\Delta t}_{h}(T)\|^{2}_{2}+\nu \int_{0}^{T}\|\nabla u^{\Delta
    t}_{h}(t)\|^{2}_{2}\,dt \leq \frac{1}{2}\|u_{0}\|^{2}_{2},
\end{equation}
a standard lower semi-continuity argument is enough to infer that the
weak solution
\begin{equation*}
  u=v=\lim_{(h,\Delta t)\to(0,0)} u^{\Delta t}_{h}=\lim_{(h,\Delta
    t)\to(0,0)} v^{\Delta t}_{h}, 
\end{equation*}
satisfies also the global energy inequality~\eqref{eq:gei}.

\medskip

\paragraph{\textbf{Step 2: Proof of the Local Energy Inequality}}

In order to conclude the proof of Theorem~\ref{teo:main} we need to
prove that the Leray-Hopf weak solution constructed in Step $1$ is
suitable. According to Definition~\ref{def:suit} this requires just to
prove the local energy inequality.  To this end, let us consider a
smooth, periodic in the space variable function $\phi\geq 0$,
vanishing for $t=0,T$; we use $P_h(u^{\Delta t}_{h}\phi)$ as test
function in the momentum equation in~\eqref{eq:thet1ac}.

The term involving the time-derivative is treated as in~\cite{BFS2019}.
  \begin{align*}
    & \int_{0}^{T}\left(\partial_t v^{\Delta t}_{h},P_h(u^{\Delta
      t}_{h}\phi)\right)\,dt
      =\int_{0}^{T}\left(\partial_t v^{\Delta t}_{h},u^{\Delta t}_{h}\phi\right)\,dt
      +\int_{0}^{T}\left(\partial_t v^{\Delta t}_{h},P_h(u^{\Delta t}_{h}\phi)-u^{\Delta
      t}_{h}\phi\right)\,dt=:I_1+I_2.
  \end{align*}
  Concerning the term $I_1$ we have that
  \begin{equation*}
    \begin{aligned}
      &=\int_0^T(\partial_t v^{\Delta t}_{h},v^{\Delta
        t}_{h})\,\phi\,dt+\int_0^T(\partial_t v^N_h,(u^{\Delta t}_{h}-v^{\Delta
        t}_{h})\,\phi)\,dt
      =:I_{11}+I_{12}.
    \end{aligned}
  \end{equation*}
  Let us first consider $I_{11}$. By splitting the integral over $[0,T]$ as the sum of
  integrals over $[t_{m-1},t_m]$ and, by integration by parts, we get
  \begin{equation*}
    \begin{aligned}
      &\int_0^T(\partial_t v^{\Delta t}_{h},v^{\Delta
        t}_{h}\phi)\,dt=\sum_{m=1}^N\int_{t_{m-1}}^{t_m}(\partial_t v^{\Delta
        t}_{h},v^{\Delta t}_{h}\phi)\,dt=
      \sum_{m=1}^N\int_{t_{m-1}}^{t_m}(\frac{1}{2}\partial_t |v^{\Delta
        t}_{h}|^2,\phi)\,dt
      \\
      &=\frac{1}{2}\sum_{m=1}^N(|u_h^{m}|^2,\phi(t_{m},x))-
      (|u_h^{m-1}|^2,\phi(t_{m-1},x))-\sum_{m=1}^N\int_{t_{m-1}}^{t_m}(\frac{1}{2}
      |v^{\Delta t}_{h}|^2,\partial_t\phi)\,dt,
    \end{aligned}
  \end{equation*}
  where we used that  $\partial_t v^{\Delta
    t}_{h}(t)=\frac{u_h^m-u_h^{m-1}}{\Delta t}$, 
  for $t\in[t_{m-1},t_m[$. Next, since the sum telescopes and $\phi$
  is with compact support in $(0,T)$ we get
  \begin{equation*}
    \int_0^T(\partial_t v^{\Delta t}_{h},v^{\Delta t}_{h}\phi)\,dt
    =-\int_{0}^{T}\big(\frac{1}{2}
    |v^{\Delta t}_{h}|^2,\partial_t\phi\big)\,dt. 
  \end{equation*}
  By the strong convergence of $v^{\Delta t}_{h}\rightarrow u$ in
  $L^2(0,T;L^2_{\#})$ we 
  can conclude that
  \begin{equation*}
    \lim_{(\Delta t,h)\to(0,0)}\int_0^T(\partial_t v^{\Delta t}_{h},v^{\Delta
      t}_{h}\phi)\,dt= -\int_{0}^{T}\big(\frac{1}{2}|u|^2,\partial_t\phi\big)\,dt.
  \end{equation*} 
  Then, we consider the term $I_{12}$. Since $u^{\Delta t}_{h}$ is constant on the
  interval $[t_{m-1}, t_m[$ we can write
  \begin{equation*}
    \begin{aligned}
      &\int_0^T(\partial_t v^{\Delta t}_{h},(u^{\Delta t}_{h}-v^{\Delta
        t}_{h})\,\phi)\,dt= -\sum_{m=1}^N\int_{t_{m-1}}^{t_m}(\partial_t(v^{\Delta
        t}_{h}-u^{\Delta t}_{h}),(v^{\Delta t}_{h}-u^{\Delta t}_{h})\,\phi)\,dt
      \\
       &=\sum_{m=1}^N\int_{t_{m-1}}^{t_m}\left(\frac{|v^{\Delta t}_{h}-u^{\Delta
            t}_{h}|^2}{2}, \partial_t\phi\right)\,dt,
 \end{aligned}
  \end{equation*}
  since the sum telescopes. Hence, we have that
  $u^{\Delta t}_{h}-v^{\Delta t}_{h}$ vanishes (strongly) in
  $L^2(0,T;L^{2}_{\#})$, provided that $\Delta t=o(h^{1/2})$. Then,
  $I_{12}\rightarrow0$ as $(\Delta t,h)\rightarrow(0,0)$.
\begin{remark}
  It is at this point that the coupling between $h$ and $\Delta t$
  plays a role. For the convergence of the other terms the discrete
  commutation property is needed. This is the reason we are skipping
  some details from the other proofs, since they are very close to
  that in the cited references. 
\end{remark}

  We have that the $I_2\to 0$ as $(\Delta t,h)\rightarrow(0,0)$. Indeed, by the
  discrete commutator property~\eqref{eq:dcp1}, Proposition~\ref{prop:1}, and the inverse
  inequality~\eqref{eq:imm} we can infer
  \begin{align*}
    \big|I_2\big|&
    \leq \int_{0}^{T}\|\partial_t v^{\Delta t}_{h}\|_{H^{-1}}\|P_h(u^{\Delta
    t}_{h}\phi)-u^{\Delta t}_{h}\phi\|_{H^1}dt
    \\
                 & \leq ch^{\frac{1}{2}}\|\partial_t v^{\Delta
                   t}_{h}\|_{L^{\frac{4}{3}}(H^{-1})}\|u^{\Delta
                   t}_{h}\|^{\frac{1}{2}}_{L^\infty(L^2)}\|u^{\Delta 
                   t}_{h}\|^{\frac{1}{2}}_{L^2(H^1)}\leq c\,h^{\frac{1}{2}}.
  \end{align*}
  Hence, also this term vanishes as $h\to0$, ending the analysis of
  the term involving the time-derivative. \par
  Concerning the viscous term, we write
  \begin{equation*}
    \begin{aligned}
      (\nabla u^{\Delta t}_{h},\nabla P_h(u^{\Delta t}_{h} \phi))
      & =(|\nabla u^{\Delta t}_{h}|^2,\phi)-(\frac{1}{2}|u^{\Delta
        t}_{h}|^2,\Delta\phi)+R_{visc},
    \end{aligned}
  \end{equation*}
  with the ``viscous remainder'' $R_{visc}:=\big(\nabla u^{\Delta
    t}_h,\nabla [P_h(u^{\Delta t}_h\phi)-u^{\Delta t}_h\phi]\big)$. 
  Since $u^{\Delta t}_{h}$ converges to $u$ weakly in
  $L^2(0,T;H_{\#}^1)$ and strongly in $L^2(0,T;L^2_{\#})$,
  \begin{equation*}
    \begin{aligned}
      \liminf_{(\Delta t,h)\to (0,0)}\int_0^T (|\nabla u^{\Delta
        t}_{h}|^2,\phi)\,dt &\geq \int_0^T (|\nabla u|^2,\phi)\, dt,
      \\
\frac{1}{2}      \int_0^T (| u^{\Delta t}_{h}|^2,\Delta\phi)\,dt &\to
\frac{1}{2}      \int_0^T (| u|^2,\Delta\phi)\,dt.
    \end{aligned}
  \end{equation*}
  For the remainder $R_{visc}$, by using again the discrete commutator
  property from Definition~\ref{def:dcp}, we have that
  \begin{equation*}
   \left| \int_{0}^{T}  R_{visc}\,dt\right|\leq c\,
   h\int_{0}^{T}\|\nabla u^{\Delta t}_{h}\|^{2}_{2}\,dt\to0 
    \qquad \textrm{ as }(\Delta t,h)\to(0,0). 
  \end{equation*}

  We consider now the nonlinear term $b_h$. We have
  \begin{equation}
    \label{eq:nonlinear_proof}
    \begin{aligned}
      b_h(u^{\Delta t}_{h},u^{\Delta t}_{h},P_h(u^{\Delta t}_{h}\phi))
      & =b_h(u^{\Delta t}_{h},u^{\Delta t}_{h},u^{\Delta t}_{h}\phi)+R_{nl}.
    \end{aligned}
  \end{equation}
  The ``nonlinear remainder'' $R_{nl}:=b_h(u^{\Delta t}_{h},u^{\Delta
    t}_{h},P_h(u^{\Delta 
    t}_{h}\phi)-u^{\Delta t}_{h}\phi)$ can be estimated by using the
  discrete commutator property,~\eqref{eq:imm}, and
 ~\eqref{eq:estimate-nonlinear1},~\eqref{eq:estimate-nonlinear1bis},
 ~\eqref{eq:estimate-nonlinear1tris} for the choices of the nonlinear
  term approximation in Case~1, Case~2, and Case~3,
  respectively. Indeed, we have 
  \begin{equation}
    \begin{aligned}
      |R_{nl}| & \leq \|nl_h(u^{\Delta t}_{h},u^{\Delta t}_{h})\|_{H^{-1}}\|P_h(u^{\Delta
        t}_{h}\phi)-u^{\Delta t}_{h}\phi\|_{H^1}
      \leq c\,\sqrt{h} \|u^{\Delta t}_{h}\|_{2}\,\|u^{\Delta t}_{h}\|_{H^1}^{2},
    \end{aligned}
  \end{equation}
  hence, by integrating in time
  \begin{equation*}
    \int_{0}^{T}R_{nl}\,dt\to0 \quad\textrm{ as }(\Delta t,h)\to(0,0).
  \end{equation*}
  The last term we consider is that involving the pressure. By
  integrating by parts we have
  \begin{equation*}
    \begin{aligned}
      (p^{\Delta t}_{h},\dive P_h(u^{\Delta t}_{h}\phi))&
      & =(p_h^{\Delta t} u^{\Delta t}_{h},\nabla\phi)+R_{p1}+R_{p2}.
    \end{aligned}
  \end{equation*}
  where the two ``pressure remainders'' are defined as follows
  \begin{equation*}
    R_{p1}:=\big(p^{\Delta t}_h,\dive(P_h(u^{\Delta t}_h\phi)-u^{\Delta
      t}_h\phi)\big)\qquad\text{and}\qquad 
    R_{p2}:=\big(\phi \, p^{\Delta t}_h,\dive u^{\Delta t}_h\big).
  \end{equation*}
By using again the discrete commutator property~\eqref{eq:dcp2}
  and~\eqref{eq:imm} we easily get
  \begin{equation*}
    \begin{aligned}
      |R_{p1}|
      &\leq c\, h\|p^{\Delta t}_{h}\|_{2}\,\|u^{\Delta t}_{h}\|_{H^1}
      \end{aligned}
    \end{equation*}
 which implies
  \begin{equation*}
    \int_0^T R_{p1}\,dt\to0\qquad  \textrm{ as
    }(\Delta t,h)\to(0,0).
  \end{equation*}
  The term $R_{p2}$ can be treated in the same way but now using the
  discrete commutation property for the projector over $Q_{h}$
  \begin{equation*}
    \begin{aligned}
      |R_{p2}|&\leq c\|Q_h(p^{\Delta t}_{h}\phi)-p^{\Delta
        t}_{h}\phi\|_{2}\,\|u^{\Delta t}_{h}\phi\|_{H^1} \leq c\,
      h^{\frac12}\|p^{\Delta t}_{h}\|_{L^{\frac43}(L^2)}\|u^{\Delta
        t}_{h}\|_{L^2(H^1)}^\frac{1}{2}\|u^{\Delta
        t}_{h}\|_{L^\infty(L^2)}^\frac{1}{2},
    \end{aligned}
  \end{equation*}
  and finally this implies that
  \begin{equation*}
    \int_0^T R_{p2}\,dt\to0\quad \textrm{ as  }(\Delta t,h)\to(0,0). 
  \end{equation*}
 The convergence  
 \begin{equation*}
 \begin{aligned}
 \int_{0}^{T}(p_h^{\Delta t} u^{\Delta
   t}_{h},\nabla\phi)\to\int_{0}^{T}(p\,  u,\nabla\phi),
 \end{aligned}
\end{equation*}
is an easy consequence of~\eqref{eq:convuv}, \eqref{eq:convp} and
Proposition~\ref{prop:1} b).  This steps are common to the three
cases.

 \bigskip
 
 We now treat the inertial term, in the Case~1 the definition of
 $nl_{h}$ in~\eqref{eq:nonlinear} allows us to handle the first term
 on the right hand side in~\eqref{eq:nonlinear_proof} with some
 integration by parts as follows:
  \begin{align*}
    b_h(u^{\Delta t}_{h},u^{\Delta t}_{h},u^{\Delta t}_{h}\phi)
     & =
       -\left(u^{\Delta t}_{h} \frac12|u^{\Delta
        t}_{h}|^2,\nabla\phi\right).
  \end{align*}
  By arguing as in~\cite{BFS2019} it can be proved that
  \begin{equation*}
    u^{\Delta t}_{h} \frac12|u^{\Delta t}_{h}|^2\to u
    \frac12|u|^2\quad \textrm{strongly in } 
    L^{1}(0,T;L^1), \quad\textrm{ as }(\Delta t,h)\to(0,0),
  \end{equation*}
  and one shows that
  \begin{equation*}
    \int_0^T  b_h(u^{\Delta t}_{h},u^{\Delta t}_{h},u^{\Delta
      t}_{h}\phi)\,dt \to -\int_0^T\Big(u
    \frac12|u|^2,\nabla\phi\Big)\,dt\quad\textrm{ as 
    }\quad(\Delta t,h)\to(0,0). 
  \end{equation*}

\bigskip

In the Case~2 the result is much simpler since, by direct computations
one shows that for smooth enough $w$ we have (by a point-wise equality,
where $\epsilon_{ijk}$ the totally anti-symmetric tensor)
\begin{equation*}
  \begin{aligned}
    \big[ ( \nabla\times w)\times w\big]\cdot(\phi\,
    w)&=\sum_{i,j,k,l}(\epsilon_{jki}-\epsilon_{jlm})
    \partial_{l}w_{m}w_{k}w_{i} \phi
    \\
    &=\phi\sum_{i,k}w_{k}\partial_{k}w_{i}w_{i}-w_{i}\partial_{i}w_{k}w_{k}=0.
  \end{aligned}
\end{equation*}
Hence, we get
\begin{equation*}
      b_h(u^{\Delta t}_{h},u^{\Delta t}_{h},u^{\Delta t}_{h}\phi)=0,
\end{equation*}
and there are not terms to be estimated. 

\bigskip

In the Case~3 we get instead (cf.~\cite[Lemma~4.1]{Gue2006})
\begin{equation*}
  \begin{aligned}
      b_h(u^{\Delta t}_{h},u^{\Delta t}_{h},u^{\Delta
        t}_{h}\phi)&=-\frac{1}{2}\big(\mathcal{K}_{h}(|u^{\Delta
        t}_{h}|^{2}),\dive(\phi\,  u^{\Delta t}_{h})\big)
      \\
      &=-\frac{1}{2}\big(u^{\Delta t}_{h}|u^{\Delta
        t}_{h}|^{2},\nabla \phi\big)+R_{1}+R_{2}
\end{aligned}
\end{equation*}
with
\begin{equation*}
  R_{1}:=-\frac{1}{2}\big(u^{\Delta t}_{h}\mathcal{K}_{h}(|u^{\Delta
    t}_{h}|^{2})-u^{\Delta t}_{h}|u^{\Delta t}_{h}|^{2},\nabla\phi\big)
  \qquad\text{and}\qquad
    R_{2}:=-\frac{1}{2}\big(\phi\,\mathcal{K}_{h}(|u^{\Delta
    t}_{h}|^{2}),\dive u^{\Delta
    t}_{h}\big),
\end{equation*}
The strong $L^{s'}(0,T;L^2)$-convergence of $\mathcal{K}_{h}(|u^{\Delta
    t}_{h}|^{2}$ implies that $\int_{0}^{T}|R_{1}|\,dt\to0$.
While using the discrete commutator property for $R_{2}$ we estimate
\begin{equation*}
  \begin{aligned}
    |R_{2}|&= \frac{1}{2} \big(\phi\,\mathcal{K}_{h}(|u^{\Delta
      t}_{h}|^{2})-Q_{h}(\phi\,\mathcal{K}_{h}(|u^{\Delta
      t}_{h}|^{2})) ,\dive u^{\Delta t}_{h}\big)|
    \\
    &\leq c h\|\mathcal{K}_{h}(|u^{\Delta t}_{h}|^{2})\|_{2}\,
    \|u^{\Delta t}_{h}\|_{H^{1}}\leq c h\|u^{\Delta t}_{h}\|_{4}^{2}\,
    \|u^{\Delta t}_{h}\|_{H^{1}}
    \\
    &\leq c h^{1/2}\|u^{\Delta t}_{h}\|_{2}^{2} \,\|u^{\Delta
      t}_{h}\|_{H^{1}}^{2},
\end{aligned}
\end{equation*}
which shows that $\int_{0}^{T}|R_{2}|dt\to0$.
\end{proof}
\section{Extension to other second order schemes}
\label{sec:6}
The techniques developed in the previous sections are general enough
to be used to handle with minor changes, also some more general second
order schemes, as for instance the Crank-Nicolson with Linear
Extrapolation~\eqref{eq:disintro3} and the
Crank-Nicholson/Adams-Bashforth~\eqref{eq:disintro2}, as reported
below. We have the following result
\begin{theorem}
  Let the same assumptions of Theorem~\ref{teo:main} be satisfied and
  replace assumption~\eqref{eq:onehalf}, by~\eqref{eq:new-restriction}
  for the~\eqref{eq:disintro3} algorithm and replace
  assumption~\eqref{eq:onehalf}, by~\eqref{eq:new-restriction2} for
  the~\eqref{eq:disintro2} algorithm.  Then, solutions of both schemes
  converge to a suitable weak solution of the NSE.
\end{theorem}

The proofs are in the same spirit of those of the previous section,
once appropriate estimates (independent of $m$) are proved. For this
reason we just include the changes with respect to the proofs of the
other cases previously treated.

\bigskip

\noindent\textbf{Crank-Nicolson with Linear Extrapolation (CNLE)}
Another scheme which is similar to the
Crank-Nicolson~\eqref{eq:disintro} in terms of theory, but better
performing in terms of numerical properties is the Crank-Nicolson with
Linear Extrapolation as introduced in Baker~\cite{Bak1973} and studied
by Ingram~\cite{Ing2013} especially in the context of non-homogeneous
Dirichlet problems.

In this case the scheme is defined, for $m\geq2$ by
\begin{equation}
  \label{eq:disintro3}\tag{CNLE}
  \begin{aligned}
    \left(d_t u^m_h,v_h\right)+\nu(\nabla u_{h}^{m,1/2},\nabla v_h)+
   \frac{1}{2} b_h(3u_{h}^{m-1}-u_{h}^{m-2},u_{h}^{m+1/2},v_h)-(
    p_{h,}^{m},\dive v_h)=0,
    \\
    ( \dive u_{h}^{m},q_h)=0,
  \end{aligned}
\end{equation}
where the operator $b_{h}(\cdot,\cdot,\cdot)$ is the same as in
``Case~1'' of the previous section and for $m=1$ the scheme is
replaced by~\eqref{eq:disintro} to be consistent with second order
time-discretization.  Scheme~\eqref{eq:disintro3} is linearly
implicit, unconditionally and nonlinearly stable, and second order
accurate, see~\cite{Bak1973,Ing2013,Ing2013b}. In~\cite{Ing2013b}, it
is shown that no time-step restriction is required for the convergence
(but with mild assumptions on the pressure) and additionally it is
proved the optimal convergence for smoother solutions.

Here we prove the following result which is not assuming any
extra-assumption neither on the Leray-Hopf weak solution $u$ nor on
the pressure $p$, and that can be used to prove the local energy
inequality, reasoning as in the previous sections.
\begin{lemma}
  \label{lem:discene2}
  Let $N\in\N$ and $m=1,..,N$. Then, for~\eqref{eq:disintro3} the
  following discrete energy-type equality holds true:
  \begin{align}
    \label{eq:abis2}
    & \frac{1}{2}(\|u_h^m\|_{2}^2-\|u_h^{m-1}\|_{2}^2)+\nu\Delta t\|\nabla
    u_{h}^{m,\frac{1}{2}}\|_{2}^2=0.
  \end{align}
  Moreover, if $u_{0}\in H^{1}_{\#}$ there exists $C>0$  such that if
  \begin{equation}
    \label{eq:new-restriction}
    \Delta t\leq\frac{\nu}{16}\min\left\{ h^{2},\frac{
        h^{3}\|u_{0}\|^{2}}{4 C^{2}}\right\},
\end{equation}
 then 
  \begin{equation*}
    \sum_{m=2}^{N}    \|u_h^m-u_h^{m-1}\|_{2}^2\leq C.
  \end{equation*}
\end{lemma}
\begin{proof}
  The first part of Lemma~\ref{lem:discene2} can be proved in a direct
  way simply using $u_{h}^{m+1/2}$ as test function,
  obtaining~\eqref{eq:abis2}; the proof of the second estimate
  requires some additional work, in the spirit
  of~\cite[Sec.~19]{MT1998}. To this end let us define
\begin{equation*}
  \delta^{m}_{h}:=\frac{u^{m}_{h}-u^{m-1}_{h}}{2},
\end{equation*}
and using $2 \Delta t\,u^{m}_{h}$ as test function
in~\eqref{eq:disintro3} we get
\begin{equation*}
  \begin{aligned}
    \|u^{m}_{h}\|^{2}_{2}-&\|u^{m-1}_{h}\|^{2}_{2}+\|\delta^{m}_{h}\|+2\nu\Delta
    t \|\nabla u^{m}_{h}\|^{2}_{2}
    \\
    & =-2\nu\Delta t(\nabla \delta^{m}_{h},\nabla u^{m})-\Delta t \,
    b_h(3u_{h}^{m-1}-u_{h}^{m-2},u_{h}^{m+1/2},u^{m}_{h}),
    \\
    & =-2\nu\Delta t (\nabla \delta^{m}_{h},\nabla u^{m})- \Delta
    t\, b_h(3
    u_{h}^{m-1}-u_{h}^{m-2},\frac{-u_{h}^{m}+u^{m-1}_{h}}{2},u^{m}_{h}),
  \end{aligned}
\end{equation*}
Hence, the right hand side can be estimated as follows
\begin{equation*}
  \begin{aligned}
    & |2\nu\Delta t(\nabla \delta^{m}_{h},\nabla u^{m})+ \Delta t\,
    b_h(3 \nabla u_{h}^{m-1}-\nabla
    u_{h}^{m-2},\delta^{m}_{h},u^{m}_{h})|,
    \\
    &\leq 2\nu\Delta t \|\nabla \delta^{m}_{h}\|_{2} \|\nabla u^{m}\|_{2}+C
    \Delta t (\|3\nabla u_{h}^{m-1}\|_{2}+\|\nabla u_{h}^{m-2}\|_{2}) \|\nabla
    u^{m}_{h}\|_{2}\| \delta^{m}_{h}\|^{1/2} _{2}\| \nabla    \delta^{m}_{h}\|^{1/2}_{2}
    \\
    &\leq \frac{2\Delta t}{h} \| \delta^{m}_{h}\|_{2}\|\nabla
    u^{m}\|_{2}+\frac{C \Delta t}{h^{3/2}} (\|3 u_{h}^{m-1}\|_{2}+\|
    u_{h}^{m-2}\|_{2}) \|\nabla u^{m}_{h}\|_{2}\| \delta^{m}_{h}\|_{2}
    \\
    & \leq \frac{1}{4}\| \delta^{m}_{h}\|^{2}_{2}+\frac{8 (\Delta
      t)^{2}}{h^{2}}\|\nabla u^{m}\|^{2}_{2}+ \frac{1}{4}\|
    \delta^{m}_{h}\|^{2}_{2}+ \frac{8C^{2} (\Delta t)^{2}}{h^{3}} (\|3
    u_{h}^{m-2}\|_{2}+ \| u_{h}^{m-2}\|_{2})^{2}\|\nabla u^{m}_{h}\|^{2}_{2} .
  \end{aligned}
\end{equation*}
Next, using the uniform estimate on $\|u_{h}^{m}\|_{2}$ coming from
the previous step we get
\begin{equation*}
  \|u^{m}_{h}\|^{2}_{2}-\|u^{m-1}_{h}\|^{2}_{2}+\frac{1}{2}\|\delta^{m}_{h}\|^{2}_{2} 
  +\nu \,   \Delta   t \|\nabla u_{h}^{m}\|^{2}_{2}
  \left(2-\frac{8\Delta t}{\nu\, h^{2}}-\frac{32 C^{2} \Delta
      t}{\nu\,h^{3}}
    \|u_{0}\|^{2}_{2}\right)\leq0,
\end{equation*}
and under the restriction on $\Delta t$ and $h$
from~\eqref{eq:new-restriction} we obtain
\begin{equation*}
  \|u^{m}_{h}\|^{2}_{2}-\|u^{m-1}_{h}\|^{2}_{2}+
  \frac{1}{2}\|\delta^{m}_{h}\|^{2}_{2}+\nu \Delta
  t \|\nabla u_{h}^{m}\|^{2}_{2}\leq0,
\end{equation*}
  which ends the proof by summation over $m$.
\end{proof}
The convergence to a weak solution satisfying the global and the local
energy inequality follows in the same manner as in~\cite{BFS2019} and
by using also the results from the previous section. Once the
estimated are proven one has just to rewrite word-by-word the proof in
the Case~1.

\bigskip

\noindent\textbf{Crank-Nicolson/Adams Bashforth} In the same spirit of
the ``{Case~1}'' we can also consider the Crank-Nicolson scheme for
the linear part and the Adams-Bashforth for the inertial one, as it is
studied for instance in~\cite[Sec.~19]{MT1998}. The algorithm reads as
follows: solve for $m\geq2$
\begin{equation}
  \label{eq:disintro2}\tag{CNAB}
  \begin{aligned}
    \left(d_t u^m_h,v_h\right)+\nu (\nabla u_{h}^{m,1/2},\nabla
    v_h)+&\frac{3}{2}b_h(u_{h}^{m-1}, u_{h}^{m-1},v_h)
    \\
    &-\frac{1}{2}b_h(u_{h}^{m-2}, u_{h}^{m-2},v_h)-&( p_{h,}^{m},\dive
    v_h)=0,&
    \\
    & &( \dive u_{h}^{m},q_h)=0,&
  \end{aligned}
\end{equation}
where $b_{h}(\cdot,\cdot,\cdot)$ is defined by means
of~\eqref{eq:nonlinear}-\eqref{eq:trilinear}, while again $u^{1}_{h}$
is obtained by an iteration of~\eqref{eq:disintro}.This method is
explicit in the nonlinear term and only conditionally
stable~\cite{MT1998,HS2007}. The~\eqref{eq:disintro2} method is
popular for approximating Navier-Stokes flows because it is fast and
easy to implement. For example, it is used to model turbulent flows
induced by wind turbine motion, turbulent flows transporting
particles, and reacting flows in complex geometries, see
Ingram~\cite{Ing2013b}.

First observe that it is possible to prove,  with a direct
argument, a sort of energy balance for the scheme,
namely an inequality of this kind
\begin{equation*}
 \frac{1}{2} \|u^{1}_{h}\|^{2}_{2}+\nu\Delta t\|\nabla
 u^{1/2}_{h}\|^{2}_{2}\leq   \frac{1}{2}\|u_{0}\|^{2}_{2}, 
\end{equation*}
but nevertheless, from the above estimate, one can also obtain by
means of the inverse inequality
\begin{equation}
  \label{eq:estimate-1}
  \frac{1}{2} \|u^{1}_{h}\|^{2}_{2}+\frac{\nu\Delta t}{4} \|\nabla
 u^{1}_{h}\|^{2}_{2}\leq  \frac{1}{2}\|u_{0}\|^{2}_{2}+\frac{\nu\Delta t}{4} \|\nabla
 u^{0}_{h}\|^{2}_{2} \leq \big(\frac{1}{2}+\frac{\nu\Delta t}{4
   h^{2}}\big) \|u_{0}\|^{2}_{2}:=K_{3}, 
\end{equation}
with $K_{3}$ independent of $\Delta t$ and $h$. The estimate for $m>1$
are obtained by a induction argument in~\cite[Lemma~19.1]{MT1998}. The
proved result is the following
\begin{lemma}
  Assume that $u_{0}\in L^{2}_{\dive}$ and~\eqref{eq:estimate-1}
  holds. Then, there exists $K_{4}$ independent of $\Delta t$ and $h$
  such that if
  \begin{equation}
    \label{eq:new-restriction2}
    \Delta t\leq \frac{4
      c_{1}^{2}}{\nu}\qquad \text{and}\qquad \frac{\Delta 
      t}{h^{3}}\leq\max\left\{\frac{1}{32\nu},\frac{c\nu}{K_{4}}\right\}, 
  \end{equation}
  then
  \begin{equation}
    \label{eq:MT}
    \begin{aligned}
      &      \|u^{n}_{h}\|^{2}_{2}\leq K_{4},
      \\
      &   \sum_{m=1}^{N}\|u^{m}_{h}-u^{m-1}_{h}\|^{2}_{2}\leq 32 K_{4},
      \\
      &  \Delta t\sum_{m=1}^{N}\|\nabla u^{m}_{h}\|^{2}_{2}\leq 4 K_{4}.
    \end{aligned}
  \end{equation}
\end{lemma}
We just comment that the proof is obtained by showing (with the same
estimates employed in the previous case) that
\begin{equation*}
  \left(1+\frac{\Delta t}{2 c_{1}^{2}}\right)\xi^{n}\leq
  \xi^{n-1}\qquad \text{where}\qquad 
  \xi^{m}:=\|u^{m}_{h}\|^{2}_{2}+\frac{1}{4}\|u^{m}_{h}-u^{m-1}_{h}\|^{2}_{2},
\end{equation*}
and then applying and inductive argument.  This is enough to prove the
standard result $u^{m}_{h}\in l^{\infty}(L^{2})\cap l^{2}(H^{1})$ from
which one deduces the estimates also on the pressure. Next passage to
the limit is again standard showing that the linear interpolated
sequence converges to a distributional solution of the NSE.

A non trivial point is to justify the global energy inequality,
because in this case the
estimate~\eqref{eq:energy-inequality-approximated} does not hold. The
functions $v^{\Delta t}_{h}$ and $u^{\Delta t}_{h}$ have the requested
regularity but do not satisfy the correct energy balance, since
\begin{equation*}
\frac{3}{2}b_h(u_{h}^{m-1},
u_{h}^{m-1}, u^{m}_{h})-\frac{1}{2}b_h(u_{h}^{m-2},
u_{h}^{m-2},u^{m}_{h})\not=0.
\end{equation*}
The correct energy balance is satisfied only in the limit
$(h,\Delta t)\to(0,0)$, but this cannot be deduced at this stage. As
usual the global energy inequality cannot be proved by means of
testing with the solution itself, but only after a limiting process,
cf.~\cite{Ber2021}.

The way of obtaining it passes through the verification that $(u,p)$
is a suitable weak solution. The validity of the local energy
inequality can be done as in~\cite{BFS2019} and results in Case~1,
once the (conditional) estimate in $\eqref{eq:MT}_2$ are
obtained. Note that in this case the restriction on the relative size
of $\Delta t$ and $h$ are needed already for the first a priori
estimate.

Next, by adapting a well-known argument in~\cite[Section~2C]{CKN1982}
we can deduce it from~\eqref{eq:lei}. In fact, it is enough to replace
$\phi$ by the product of $\phi$ and $\chi_{\epsilon}$ (which is a
mollification of $\chi_{[t_{1},t_{2}]}(t)$, the characteristic
function of $[t_{1},t_{2}]$) and pass to the limit as $\epsilon\to0$
to get
\begin{equation*}
\begin{aligned}
  &\int_{\tore}|u(t_{2})|^{2}\phi(t_{2})\,dx+ \nu
  \int_0^T\int_{\tore}|\nabla u|^2\phi\,dx dt\leq
  \\
  &\quad \leq \int_{\tore}|u(t_{1})|^{2}\phi(t_{1})\,dx +
  \int_{0}^{T}\int_{\tore}\left[\frac{|u|^
      2}{2}\left(\partial_t\phi+\nu\Delta\phi\right)%
    +\left(\frac{|u|^2}{2}+p\right)u\cdot\nabla \phi\right]\,dxdt,
\end{aligned}
\end{equation*}
and the above formula is particularly significant if
$\phi(\tau,x)\not =0$ in $(t_{1},t_{2})$. Next in the above formula
one can take a sequence $\phi_{n}$ of smooth functions converging to
the function $\phi\equiv1$ and at least in the whole space or in the
space periodic setting one gets the global energy
inequality~\eqref{eq:gei} as is explained at the beginning
of~\cite[Section~8]{CKN1982}. Moreover, the same argument applied to
arbitrary time intervals shows also that
  \begin{equation*}
    \frac12\|u(t_{2})\|_{2}^2+\nu\int_{t_{1}}^{t_{2}}\|\nabla
    u(s)\|_{2}^2\,ds\leq\frac12 
    \|u(t_{1})\|_{2}^2\qquad\textrm{ for all } 0\leq t_{1}\leq
    t_{2}\leq T. 
  \end{equation*}
hence that the strong global energy inequality holds true.

\begin{thebibliography}{10}

\bibitem{ABC2021} {\sc D.~Albritton, E.~Bru\'e, and M.~Colombo}, {\em Non-uniqueness of Leray solutions of the forced Navier-Stokes equations}, Arxiv 2112:03116, (2021).

\bibitem{Bak1973}
{\sc G.~A. Baker}, {\em P{rojection} {methods} {for} {boundary}-{value}
  {problems} {for} {equations} {of} {elliptic} {and} {parabolic} {type} {with}
  {discontinuous} {coefficients}}, ProQuest LLC, Ann Arbor, MI, 1973,
\newblock Thesis (Ph.D.)--Cornell University.

\bibitem{Ber2018}
{\sc L.~C. Berselli}, {\em Weak solutions constructed by semi-discretization
  are suitable: the case of slip boundary conditions}, Int. J. Numer. Anal.
  Model., 15 (2018), pp.~479--491.

\bibitem{Ber2021}
{\sc L.~C. Berselli}, {\em Three-dimensional {N}avier-{S}tokes equations for
  turbulence}, Mathematics in Science and Engineering, Academic Press, London,
  [2021] \copyright 2021.

\bibitem{BFS2019}
{\sc L.~C. Berselli, S.~Fagioli, and S.~Spirito}, {\em Suitable weak solutions
  of the {N}avier-{S}tokes equations constructed by a space-time numerical
  discretization}, J. Math. Pures Appl. (9), 125 (2019), pp.~189--208,

\bibitem{BS2012}
{\sc L.~C. Berselli and S.~Spirito}, {\em On the vanishing viscosity limit of
  3{D} {N}avier-{S}tokes equations under slip boundary conditions in general
  domains}, Comm. Math. Phys., 316 (2012), pp.~171--198,
  
\bibitem{BS2014}
{\sc L.~C. Berselli and S.~Spirito}, {\em An elementary approach to the
  inviscid limits for the 3{D} {N}avier-{S}tokes equations with slip boundary
  conditions and applications to the 3{D} {B}oussinesq equations}, NoDEA
  Nonlinear Differential Equations Appl., 21 (2014), pp.~149--166,
  
\bibitem{BS2016}
{\sc L.~C. Berselli and S.~Spirito}, {\em Weak solutions to the
  {N}avier-{S}tokes equations constructed by semi-discretization are suitable},
  in Recent advances in partial differential equations and applications,
  vol.~666 of Contemp. Math., Amer. Math. Soc., Providence, RI, 2016,
  pp.~85--97,
  
\bibitem{BX191}
{\sc J.~H. Bramble and J.~Xu}, {\em Some estimates for a weighted {$L^2$}
  projection}, Math. Comp., 56 (1991), pp.~463--476,

\bibitem{BS2008}
{\sc S.~C. Brenner and L.~R. Scott}, {\em The mathematical theory of finite
  element methods}, vol.~15 of Texts in Applied Mathematics, Springer, New
York, third~ed., 2008,

\bibitem{CKN1982}
{\sc L.~Caffarelli, R.~Kohn, and L.~Nirenberg}, {\em Partial regularity of
  suitable weak solutions of the {N}avier-{S}tokes equations}, Comm. Pure Appl.
  Math., 35 (1982), pp.~771--831.

\bibitem{Car2002}
{\sc C.~Carstensen}, {\em Merging the {B}ramble-{P}asciak-{S}teinbach and the
  {C}rouzeix-{T}hom\'{e}e criterion for {$H^1$}-stability of the
  {$L^2$}-projection onto finite element spaces}, Math. Comp., 71 (2002),
  pp.~157--163.

\bibitem{DST2021}
{\sc L.~Diening, J.~Storn, and T.~Tscherpel}, {\em On the {S}obolev and
  {$L^p$}-stability of the {$L^2$}-projection}, SIAM J. Numer. Anal., 59
  (2021), pp.~2571--2607.

\bibitem{DDW1974}
{\sc J.~Douglas, Jr., T.~Dupont, and L.~Wahlbin}, {\em The stability in
  {$L^{q}$} of the {$L^{2}$}-projection into finite element function spaces},
  Numer. Math., 23 (1974/75), pp.~193--197.

\bibitem{Gue2006}
{\sc J.-L. Guermond}, {\em Finite-element-based {F}aedo-{G}alerkin weak
  solutions to the {N}avier-{S}tokes equations in the three-dimensional torus
  are suitable}, J. Math. Pures Appl. (9), 85 (2006), pp.~451--464.

\bibitem{Gue2007}
{\sc J.-L. Guermond}, {\em Faedo-{G}alerkin weak solutions of the
  {N}avier-{S}tokes equations with {D}irichlet boundary conditions are
  suitable}, J. Math. Pures Appl. (9), 88 (2007), pp.~87--106.

\bibitem{Gue2008}
{\sc J.-L. Guermond}, {\em On the use of the notion of suitable weak solutions
  in {CFD}}, Internat. J. Numer. Methods Fluids, 57 (2008), pp.~1153--1170.

\bibitem{GOP2004}
{\sc J.-L. Guermond, J.~T. Oden, and S.~Prudhomme}, {\em Mathematical
  perspectives on large eddy simulation models for turbulent flows}, J. Math.
  Fluid Mech., 6 (2004), pp.~194--248.

\bibitem{HL1996}
{\sc Y.~He and K.~Li}, {\em Nonlinear {G}alerkin method and two-step method for
  the {N}avier-{S}tokes equations}, Numer. Methods Partial Differential
  Equations, 12 (1996), pp.~283--305.

\bibitem{HS2007}
{\sc Y.~He and W.~Sun}, {\em Stability and convergence of the
  {C}rank-{N}icolson/{A}dams-{B}ashforth scheme for the time-dependent
  {N}avier-{S}tokes equations}, SIAM J. Numer. Anal., 45 (2007), pp.~837--869.

\bibitem{HR1990}
{\sc J.~G. Heywood and R.~Rannacher}, {\em Finite-element approximation of the
  nonstationary {N}avier-{S}tokes problem. {IV}. {E}rror analysis for
  second-order time discretization}, SIAM J. Numer. Anal., 27 (1990).

\bibitem{Hor1987}
{\sc K.~Horiuti}, {\em Comparison of conservative and rotational forms in large
  eddy simulation of turbulent channel flow}, J. Comput. Phys.,
  71 (1987), pp.~343--370.

\bibitem{Ing2013}
{\sc R.~Ingram}, {\em A new linearly extrapolated {C}rank-{N}icolson
  time-stepping scheme for the {N}avier-{S}tokes equations}, Math. Comp., 82
  (2013), pp.~1953--1973.

\bibitem{Ing2013b}
{\sc R.~Ingram}, {\em Unconditional convergence of high-order extrapolations of
  the {C}rank-{N}icolson, finite element method for the {N}avier-{S}tokes
  equations}, Int. J. Numer. Anal. Model., 10 (2013), pp.~257--297.

\bibitem{LMNOR2009}
{\sc W.~Layton, C.~C. Manica, M.~Neda, M.~Olshanskii, and L.~G. Rebholz}, {\em
  On the accuracy of the rotation form in simulations of the {N}avier-{S}tokes
  equations}, J. Comput. Phys., 228 (2009), pp.~3433--3447.

\bibitem{Lio1998}
{\sc P.-L. Lions}, {\em Mathematical topics in fluid mechanics. {V}ol. 2},
  vol.~10 of Oxford Lecture Series in Mathematics and its Applications, The
  Clarendon Press, Oxford University Press, New York, 1998.
\newblock Compressible models, Oxford Science Publications.

\bibitem{MT1998}
{\sc M.~Marion and R.~Temam}, {\em Navier-{S}tokes equations: theory and
  approximation}, in Handbook of numerical analysis, {V}ol. {VI}, Handb. Numer.
  Anal., VI, North-Holland, Amsterdam, 1998, pp.~503--688.

\bibitem{QV1994}
{\sc A.~Quarteroni and A.~Valli}, {\em Numerical approximation of partial
  differential equations}, vol.~23 of Springer Series in Computational
  Mathematics, Springer-Verlag, Berlin, 1994.

\bibitem{Sche1977}
{\sc V.~Scheffer}, {\em Hausdorff measure and the {N}avier-{S}tokes equations},
  Comm. Math. Phys., 55 (1977), pp.~97--112.

\bibitem{Tem1977b}
{\sc R.~Temam}, {\em Navier-{S}tokes equations. {T}heory and numerical
  analysis}, North-Holland Publishing Co., Amsterdam-New York-Oxford, 1977.
\newblock Studies in Mathematics and its Applications, Vol. 2.

\bibitem{Tho1997}
{\sc V.~Thom{\'e}e}, {\em Galerkin finite element methods for parabolic
  problems}, vol.~25 of Springer Series in Computational Mathematics,
  Springer-Verlag, Berlin, 1997.

\bibitem{Ton2004}
{\sc F.~Tone}, {\em Error analysis for a second order scheme for the
  {N}avier-{S}tokes equations}, Appl. Numer. Math., 50 (2004), pp.~93--119.

\bibitem{Vas2007}
{\sc A.~F. Vasseur}, {\em A new proof of partial regularity of solutions to
  {N}avier-{S}tokes equations}, NoDEA Nonlinear Differential Equations Appl..

\bibitem{Zang1991} {\sc T.~Zang}, {\em On the rotation and
    skew-symmetric forms for incompressible flow simulations},
  Appl. Numer. Math., 7 (1991), pp.~27--40.

\end{thebibliography}
\section*{Acknowledgments}
The authors thank V. DeCaria for useful suggestions and comments on an
early draft of the paper. The authors acknowledge support by
INdAM-GNAMPA.

\def\ocirc#1{\ifmmode\setbox0=\hbox{$#1$}\dimen0=\ht0 \advance\dimen0
  by1pt\rlap{\hbox to\wd0{\hss\raise\dimen0
  \hbox{\hskip.2em$\scriptscriptstyle\circ$}\hss}}#1\else {\accent"17 #1}\fi}
  \def\polhk#1{\setbox0=\hbox{#1}{\ooalign{\hidewidth
  \lower1.5ex\hbox{`}\hidewidth\crcr\unhbox0}}} \def\cprime{$'$}

\end{document}